\documentstyle[amscd,amssymb,verbatim]{amsart}
\newcommand{\lrar}[1]{\begin{picture}(50,10)(-25,-5)                          
\put(-25,0){\vector(1,0){50}}
\put(0,5){\makebox(0,0)[b]{\mbox{$#1$}}}
\end{picture}}

\newcommand{\ldar}[1]{\begin{picture}(10,50)(-5,-25)
\put(0,25){\vector(0,-1){50}}
\put(5,0){\mbox{$#1$}}
\end{picture}}

\newcommand{\bff}{{\bf f}}
\newcommand{\sgn}{\operatorname{sgn}}

\newcommand{\vol}{\operatorname{vol}}

\newcommand{\Spin}{\operatorname{Spin}}
\newcommand{\SO}{\operatorname{SO}}

\newcommand{\ssl}{\operatorname{sl}}
\newcommand{\GL}{\operatorname{GL}}
\newcommand{\Cr}{\operatorname{Cr}}

\newcommand{\Sym}{\operatorname{Sym}}

\newcommand{\disj}{\bigsqcup}

\newcommand{\coker}{\operatorname{coker}}
\newcommand{\DD}{{\cal D}}

\newcommand{\OO}{{\cal O}}
\newcommand{\Gal}{\operatorname{Gal}}

\newcommand{\SL}{\operatorname{SL}}

\newcommand{\Aut}{\operatorname{Aut}}

\newcommand{\G}{{\Bbb G}}
\newcommand{\A}{{\Bbb A}}

\newcommand{\lan}{\langle}
\newcommand{\ran}{\rangle}

\newcommand{\supp}{\operatorname{supp}}

\newcommand{\Th}{\Theta}

\newcommand{\inv}{\operatorname{inv}}

\newcommand{\de}{\delta}
\newcommand{\eps}{\epsilon}
\newcommand{\We}{\bigwedge}
\renewcommand{\ker}{\operatorname{ker}}

\numberwithin{equation}{subsection}

\newtheorem{thm}{Theorem}[subsection]
\newtheorem{prop}[thm]{Proposition}
\newtheorem{lem}[thm]{Lemma}

\newenvironment{rem}{\vspace{3mm}\noindent
{\bf Remark.}}{\vspace{3mm}}

\newcommand{\Pf}{\noindent {\it Proof}}
\newcommand{\id}{\operatorname{id}}

\newcommand{\ov}{\overline}
\newcommand{\we}{\wedge}

\newcommand{\rk}{\operatorname{rk}}
\newcommand{\ra}{\rightarrow}

\newcommand{\FF}{{\cal F}}

\newcommand{\SS}{{\cal S}}

\renewcommand{\O}{{\cal O}}

\newcommand{\Hom}{\operatorname{Hom}}

\renewcommand{\a}{\alpha}
\renewcommand{\b}{\beta}
\newcommand{\om}{\omega}

\newcommand{\la}{\lambda}
\newcommand{\th}{\theta}
\newcommand{\C}{{\Bbb C}}
\newcommand{\R}{{\Bbb R}}
\newcommand{\Z}{{\Bbb Z}}
\newcommand{\Q}{{\Bbb Q}}

\newcommand{\Ga}{\Gamma}

\newcommand{\wt}{\widetilde}

\newcommand{\sign}{\operatorname{sign}}

\newcommand{\ed}{\qed\vspace{3mm}}
\newcommand{\bG}{{\bf G}}
\newcommand{\bL}{{\bf L}}
\newcommand{\bS}{{\bf S}}
\newcommand{\bI}{{\bf I}}
\newcommand{\bU}{{\bf U}}
\newcommand{\bC}{{\bf C}}
\newcommand{\bZ}{{\bf Z}}
\newcommand{\bX}{{\bf X}}
\newcommand{\bV}{{\bf V}}
\newcommand{\bW}{{\bf W}}
\newcommand{\bH}{{\bf H}}
\newcommand{\bT}{{\bf T}}
\newcommand{\bK}{{\bf K}}

\renewcommand{\gg}{{\frak g}}
\newcommand{\cc}{{\frak c}}
\newcommand{\hh}{{\frak h}}

\title{Generalization of a theorem of Waldspurger to nice representations}
\author{D. Kazhdan and A. Polishchuk}
\thanks{Both authors were supported in part by NSF grants}

\begin{document}
\maketitle

The goal of this paper is to present some examples of equalities of
integrals over local fields predicted by the
stationary phase approximation. These examples
should be considered in the context of the algebraic integration
theory proposed by the first author in \cite{Ka}.
From the point of view of this theory,
our main theorem gives examples of pairs of algebro-geometric data
over a given local field $E$ producing equal integrals over
arbitrary finite extensions of $E$.
Another motivation for this work is the theory of
Sato's functional equations associated with
prehomogeneous vector spaces over local fields. As an application of
our techniques we find an explicit form of these equations for
the action of $\GL_n$ on symmetric $n\times n$ matrices in the case
when $n$ is odd (this is a generalization of a particular case of the
equations obtained by W.~J.~Sweet Jr. in \cite{Sw}).

Let $\bG$ be a (connected) simply connected semisimple algebraic group over
a local field $E$ of characteristic zero,
$\rho:\bG\ra\Aut(\bV)$ be a
rational representation defined over $E$.

\noindent
{\bf Definition}. A representation $\rho$ is called {\it nice} if 
the generic stabilizer subgroup $\bH$ is connected and reductive.
\footnote{This notion is slightly more general than
the one defined in section 7 of \cite{Ka}.}

\vspace{2mm}

Note that the adjoint representation of $\bG$ is nice. Other
interesting examples
can be found using Elashvili's tables of representations
with positive-dimensional generic stabilizers (see \cite{El1},\cite{El2}).

Let $\bV_0\subset \bV$ be a non-empty Zariski open affine subset
such that the stabilizer of any point in $\bV_0$
is conjugate to $\bH$ over $\ov{E}$, $\bV^{\vee}_0\subset\bV^{\vee}$
be a similar subset in the dual representation.
We are interested in the subspaces $\DD^{st}$, $\DD^{st}_{\eps}$ of
$\bG(E)$-invariant
distributions on $\bV(E)$ consisting of {\it stable} and $\eps$-stable
distributions. To define these spaces we have to introduce some
notations.

Recall that for a non-degenerate quadratic form $q$ over a local field
$E$ the Hasse-Witt invariant $\eps(q)=\pm 1$ is defined as follows:
choose coordinates $x_1,\ldots,x_n$ in such a way
that $q=a_1x_1^2+\ldots+a_nx_n^2$ and set $\eps(q)=\prod_{i<j}(a_i,a_j)$
where $(\cdot,\cdot)$ denotes the Hilbert symbol (if rank of $q$ is equal
to $1$ we set $\eps(q)=1$). Note that if $E=\R$ then
$\eps(q)=(-1)^{i(i-1)/2}$ where $i$ is the number of negative squares
in $q$ (for $E=\C$ we have $\eps(q)=1$).
For a pair of non-degenerate quadratic forms $q$ and $q'$
we define the relative Hasse-Witt invariant $\eps(q,q'):=\eps(q)\eps(q')$.

Let $\gg$ be the Lie algebra of $\bG$, $Q$ be the Killing form
on $\gg$.
For every point $x\in\bV_0(E)$ we denote by $\bH_x\subset\bG$
the stabilizer of $x$, by $\hh_x\subset\gg$ its Lie algebra.
It is easy to see that the form $Q_x:=Q|_{\hh_x}$ is still
non-degenerate (see lemma \ref{Killem}).
Now for every $\bG$-orbit $\OO\subset\bV_0$ we define  
the relative sign function on $\OO(E)\times\OO(E)$ by
$$\eps(x,x'):=\eps(Q_x,Q_{x'}).$$
Sometimes we will also use the absolute sign function
$$\eps(x):=\eps(Q_x)$$
such that $\eps(x,x')=\eps(x)\eps(x')$.
Let us denote by $p:\bV_0\ra\bV_0/\bG$ the natural projection
to the geometric quotient by the action of $\bG$.

\noindent
{\bf Definition}. Let  $\delta$ be a  distribution which
is given by a locally $L^1$ function $f$:
$\delta =f|dv|$ where $|dv|$ is
the Haar measure on  $\bV(E)$. 

\noindent
a)We say that $\delta$ is {\it stable} if the function $f$ is
constant on the fibers of the natural projection 
$p(E):\bV_0(E)\ra (\bV_0/\bG)(E)$ 

\noindent
b) We say that  $\delta$ is $\eps$-stable if for (almost) any pair
$x,x'\in  \bV_0(E)$ such that $p(E)(x)=p(E)(x')$ we have 
$f(x)=\eps (x,x')f(x')$

\noindent
By the definition, the space of all stable (resp. $\eps$-stable)
distributions on $\bV(E)$ is the closure of the space of locally
$L^1$ stable (resp. $\eps$-stable) distributions
(another definition will be given
in section \ref{stabledistrsec}).

\vspace{2mm}

Our main result is that in the case when $E$ is $p$-adic, $\rho$ is nice and either
$\bG$ is simple or $\bH$ is semisimple,
the notions of stable and $\eps$-stable
distributions get switched by the Fourier transform.
Note that in the case when $\rho$ is the adjoint representation we
have $\eps\equiv 1$. Thus, in this case our result is that
stability is preserved under Fourier transform. This was proven 
previously by J.-L.~Waldspurger (see \cite{Wa1}, Cor. 1.6).
We conjecture that similar result holds in the case $E=\R$.
We also introduce the ``complementary" notions of antistable
and $\eps$-antistable distributions and prove that they get switched
by the Fourier transform in the case when $\rho$ is nice, $E$ is
$p$-adic and $\bH$ is semisimple.

The proof combines some local computations with global considerations.
\footnote{The sketch of the proof in the case $\eps\equiv 1$
was given in \cite{Ka}.}
Our main local tool is the stationary phase approximation applied
to linear functionals on $\bG$-orbits.
The main global tool we use is the stabilization techniques
introduced by R.~Kottwitz in \cite{K2}
to stabilize the elliptic semisimple part of the trace formula for $\bG$.
More precisely, we prove an analogue of the Kottwitz stabilization formula 
(Theorem 9.6 of \cite{K2}) for arbitrary nice representations over number field.

Using the same method we prove the following result concerning stable
distributions and inner forms
(sketched in \cite{Ka} in the case $\eps\equiv 1$).
Let $\bG'\ra\Aut(\bV')$ be an inner form
of a nice representation $\bG\ra\Aut(\bV)$ over $E$.
Then $\bG$-orbits on $\bV$ are in bijection with
$\bG'$-orbits on $\bV'$. For any $\bar{x} \in \bV /\bG (E)$ we denote by 
$\OO_{\bar{x}} \subset \bV ,\OO'_{\bar{x}} \subset \bV'$ the preimages of
$\bar{x}$. For any $x \in \OO_{\bar{x}} (E),x' \in \OO'_{\bar{x}}(E)$ we
define as before $\eps (x, x')=\eps(Q_x,Q'_{x'})$
(where $Q'_{x'}$ is the restriction of the Killing form $Q'$ on $\gg'$
to $\hh_{x'}$).

\noindent
{\bf Definition}. a) Let $\delta =f|dv|$, $\delta' =f'|dv'|$ be stable
distributions on $\bV(E)$, $\bV '(E)$ given by locally $L^1$ functions.
We write $ \delta \sim \delta '$ if for any $\bar x \in
\bV /\bG (E), x\in \OO_{\bar{x}} (E),x' \in \OO'_{\bar{x}} (E)$ we have
$f(x)=f'(x')$.

\noindent
b) Let $\delta =f|dv|$, $\delta' =f'|dv'|$ be $\eps$-stable
distributions on $\bV(E)$, $\bV'(E)$ given by locally $L^1$ functions. 
We write $ \delta \sim_{\eps} \delta'$ if for any $\bar x \in
\bV /\bG (E),x \in \OO_{\bar{x}} (E),x' \in \OO'_{\bar{x}} (E)$ we have
$f(x)=f'(x')\eps (x, x')$.

We extend these definitions to arbitrary stable (resp. $\eps$-stable)
distributions and prove that in the case when $E$ is $p$-adic and $\bG$ is simple,
for any stable
distributions $\delta$ and $\delta'$ on  $\bV(E)$ and $\bV'(E)$
respectively, such that $\delta \sim \delta'$, we have
$\FF(\delta) \sim_{\eps} \kappa(\bG,\bG')\FF(\delta')$ where
$\kappa(\bG,\bG')$ is the relative Hasse-Witt invariant
of Killing forms on $\gg$ and $\gg'$. The latter sign coincides with
the sign considered by Kottwitz in \cite{K0}.

Finally, we apply our main theorem to derive the following equation
for distributions on the space $\Sym_n(E)$ of symmetric $n\times n$ matrices
over a $p$-adic field $E$ in the case when $n$ is odd:
\begin{equation}\label{discreq}
\FF(\chi(\det))=c(\chi)\cdot\eps\cdot (\det,-1)^{\frac{n-1}{2}}\cdot
|\det|^{-\frac{n+1}{2}}\chi^{-1}(\det)
\end{equation}
where $\chi$ is a generic multiplicative character of a local field $E$,
$\chi(\det)$ is considered as a distribution on $\Sym_n(E)$ (defined
by analytic continuation), $\eps$ is the function on the complement
to the hypersurface $(\det=0)$ in $\Sym_n(E)$ that
assigns to a symmetric matrix the Hasse-Witt invariant of the corresponding
quadratic form, $c(\chi)$ is a non-zero constant. In the case
when $\chi=|\cdot|^s\delta$, where $\delta$ is a character of
$E^*/(E^*)^2$, the equation (\ref{discreq}) follows from Proposition 4.8 of
\cite{Sw}. We show that in the case $E=\R$ the equation (\ref{discreq}) follows from
computations of Shintani in \cite{Sh}. This confirms the conjecture that our main theorem holds also
in the real case.

The appearence of the sign $\eps$ in the above results
is due to the connection of the relative Hasse-Witt invariant
with Weil constants for quadratic forms.
First, let us introduce some more notation concerning quadratic forms.
Considering a non-degenerate quadratic form $q$
on a $E$-vector space $L$ as a symmetric isomorphism
$L\ra L^{\vee}$ we can define $\det(q)\in\det(L)^{\otimes(-2)}$,
where $\det(L)=\We^{\dim L}L$. By choosing a trivialization of $\det(L)$
we get an element of $E^*$. The corresponding element in $E^*/(E^*)^2$
does not depend on a choice of trivialization (we call it {\it determinant of
$q$ modulo squares}).
If $E$ is $p$-adic, then two non-degenerate quadratic forms
$q$ and $q'$ are equivalent over $E$ if and only if
$\rk q=\rk q'$, $\det(q)\equiv\det(q')\mod (E^*)^2$ and $\eps(q,q')=1$.
More generally, if $E\neq\C$,
then these invariants have the following $K$-theoretic meaning.
Let $W_E$ be the Witt ring of $E$. By the definition,  $W_E$ is
the abelian group generated by pairs $(L,q)$ where $L$ is a
finite-dimensional  $E$-vector space and $q$ is a nondegenerate
quadratic form on $L$ and the product in  $W_E $ comes from the
operation of tensor product  $(L',q')\cdot (L'',q''):=(L,q)$ where
$L=L'\otimes L''$. Let $W_E^1\subset W_E $ be the ideal of forms of
even dimension, $W_E^i :=(W_E ^1)^i$. We define $\ov{W}_E^i :=W_E^i
/W_E^{i+1}$. As is well-known $\ov{W}_E ^0 =\Z/2\Z$ ,$\ov{W}_E^1
=E^{*}/(E^{*})^2$, $\ov{W}_E ^2\simeq K_2(E)/2K_2(E)\simeq\Z/2\Z$.
For a pair of non-degenerate quadratic forms $(q,q')$
of the same rank, we have $(q)-(q')\in W_E^1$ and its image 
in $\ov{W}_E^1$ is identified with
$\det(q)/\det(q')\mod(E^*)^2$. If the latter element is trivial, then
we have $(q)-(q')\in W_E^2$ and its image in $\ov{W}_E ^2$
can be identified with the
relative Hasse-Witt invariant $\eps(q,q')$.
On the other hand, for a fixed (non-trivial) additive
character $\psi:E\ra U(1)$ and a non-degenerate quadratic form $q$,
A.~Weil introduced in \cite{We} a constant
$\gamma(q,\psi)$ which is a root of unity of order $8$ depending
only on the equivalence class of $q$ (and on $\psi$).
Note that the character $\psi$ defines a canonical Haar measure on $E$
which is self-dual for the Fourier transform.
Let $V$ be a vector space over $E$ on which $q$ is defined. For
a non-zero top-degree form $\nu$ on $V$ we have the induced Haar measure
$|\nu|$ on $V$ (corresponding to the self-dual Haar measure on $E$).
The Weil's constant is defined by the functional equation for distributions
\begin{equation}\label{Weileq}
\FF(\psi(q))=\gamma(q,\psi)\cdot |\det(q)/\nu^2|^{-\frac{1}{2}}
\psi(-q^{\vee}),
\end{equation}
where $q^{\vee}$ is the dual quadratic form on $V^{\vee}$,
$\psi(q)$ and $\psi(-q^{\vee})$ are considered as distributions on $V$ and 
$V^{\vee}$. By the definition,
the map $q\mapsto\gamma(q,\psi)$ extends to a homomorphism
from $W_E$ to roots of unity of order $8$. As shown in \cite{We} (nos. 25--28),
for quadratic forms $q$ and $q'$ such that $(q)-(q')\in W_E^2$
one has
$$\gamma((q)-(q'),\psi)=\eps(q,q').$$


Here is the plan of the paper. In section \ref{repthsec}
we gather some algebraic facts about nice representations.
In particular, we prove that the generic stabilizers of  
a nice representation and of its dual are conjugate to each other
(proposition \ref{dualnice}).
In section \ref{localsec} we draw consequences from the
stationary phase approximation in the case when the ground field is $p$-adic.
In section \ref{globalsec} we prove an analogue
of the stable trace formula for nice representations
over number fields 
and combine it with local information from section \ref{localsec}
to derive the main result.
Finally, in section \ref{prehom} we derive the equation (\ref{discreq}).

\noindent
{\it Notation}. All our fields have characteristic zero.
$E$ always denotes a local field, while
$F$ always denotes a number field. By $p$-adic field we mean a finite
extension of $\Q_p$. 
For a vector space $V$ over a local field $E$ and
an open subset $U\subset V$ we denote by
$\SS(U)$ the space of functions with support in $U$ which belong
to the Schwartz-Bruhat space of $V$.
When we work over a local field $E$, we fix a non-trivial additive
character $\psi:E\ra U(1)$. We fix the Haar measure on $E$ which is
self-dual with respect to the Fourier transform defined in terms of $\psi$.
For a smooth variety $\bX$ over $E$ and a non-vanishing top-degree
form $\om$ on $E$, we denote by $|\om|$ the measure on $\bX(E)$ corresponding
to $\om$ and to the above Haar measure on $E$. 
For a variety $\bX$ defined over a field $k$ and an extension
of fields $k\subset k'$, we denote by $\bX_{k'}$ the variety over $k'$
obtained from $\bX$ by the extension of scalars.
For an algebraic group $\bH$ we denote by $Z(\bH)$ its centre.
The group $\bG$ is assumed to be (connected) semisimple and 
simply connected, $\gg$ denotes the Lie algebra of $\bG$, $Q$ is
the Killing form on $\gg$.
When $\bG$ acts on a vector space $\bV$, for every point
$x\in\bV$ we denote by $\OO_x\subset\bV$ the $\bG$-orbit of $x$ and by
$\bH_x\subset\bG$ (resp. $\hh_x\subset\gg$)
the stabilizer subgroup (resp. subalgebra) of $x$.
For a field $k$ we denote by $\Ga_k$ the Galois group $\Gal(\ov{k}/k)$.
For an algebraic group $\bH$ defined over $k$ we set
$H^i(k,\bH)=H^i(\Ga_k,\bH(\ov{k}))$
(where $i\le 1$ if $\bH$ is noncommutative).

{\it Acknowledgment}. We would like to thank M.~Borovoi,
A.~Elashvili, P.~Etingof, V.~Kac and G.~Seitz for helpful discussions.

\section{Algebraic results}\label{repthsec}

Throughout this section $\bG$ is a simply connected semisimple group over
a field $k$ of characteristic zero,
$\rho:\bG\ra\Aut(\bV)$ is a
rational representation defined over $k$. In this section we
gather results we need that can be proven after passing to an
algebraic closure of $k$.

\subsection{Representations with reductive generic stabilizer}

The first restriction we impose in order for a representation
$\rho$ to be nice is the reductivity of $\hh$.
In this subsection we discuss some consequences of this condition.

It is well-known (see \cite{Lu}, \cite{R})
that there exists a non-empty $\bG$-invariant
Zariski open subset $\bV_0\subset \bV$
and a subgroup $\bH\subset\bG$ such that for every point
$x\in\bV_0(\ov{k})$ the stabilizer subgroup $\bH_x$ of $x$
is conjugate to $\bH$. Furthermore, clearly we can assume that
$\bV_0$ is invariant under the natural $\bG_m$-action on $\bV$.
The following lemma can be found in \cite{El1}
(see also \cite{PV}, Theorem 7.3).

\begin{lem}\label{sumlem}
For every $x\in\bV_0$ one has
$$\bV=\bV^{\hh_x}+\gg x.$$
\end{lem}

\Pf . This follows immediately from the surjectivity
of the map $\bG/\bH_x\times(\bV_0^{\hh_x})\ra\bV_0$,
where $\bV_0^{\hh_x}=\bV^{\hh_x}\cap\bV_0$.
\ed

Let $\hh\subset\gg$ be
the Lie algebra of $\bH$.
We will call $\bH$ (resp. $\hh$) the {\it generic
stabilizer subgroup} (resp. {\it subalgebra}) for $\rho$.
Henceforth, we always assume that $\hh$ is reductive.
The theorem of V.~L.~Popov \cite{Popov}
states that in this case generic $\bG$-orbits in $\bV$ are closed. 

\begin{prop}\label{dual} 
Let $\rho^{\vee}:\bG\ra\Aut(\bV^{\vee})$
be the dual representation to $\rho$. Then the generic
stabilizer subalgebra for $\rho^{\vee}$ is conjugate to $\hh$
over $\ov{k}$.
\end{prop}

\Pf . Clearly, we can assume that $k=\C$. Let 
$C\subset \bG(\C)$ be a maximal compact subgroup. It is well-known
that there exists a $C$-invariant positive-definite Hermitian
form $H$ on $\bV$. Let $\OO\subset\bV$ be a generic orbit.
Let us consider the restriction of the function
$x\mapsto H(x,x)$ to $\OO(\C)$.
Since $\OO$ is closed, there exists a vector $x\in\OO(\C)$ minimizing
this function. In particular, we have $H(\gg x,x)=0$.
Let us consider the functional $x^{\vee}\in\bV^{\vee}(\C)$ given by
$x^{\vee}=H(?,x)$. We claim that the stabilizer subalgebra of
$x^{\vee}$ coincides with $\hh_x$. Indeed, first we note that
$x^{\vee}\in (\gg x)^{\perp}$. According to lemma \ref{sumlem},
$\hh_x$ acts trivially on $V/\gg x$. Hence, it also
acts trivially $(\gg x)^{\perp}$, so $\hh_x x^{\vee}=0$. Let
us define the $\C$-bilinear form on $\gg/\hh_x$ by setting
$$B(\xi_1,\xi_2)=\lan \xi_1 x^{\vee},\xi_2 x\ran=-H(\xi_1\xi_2 x,x)$$
where $\xi_1,\xi_2\in\gg/\hh_x$.
Note that since $H([\xi_1,\xi_2]x,x)=0$, the form $B$ is symmetric.
Now let $\cc\subset\gg$ be the Lie algebra of $C$. We claim that
the restriction of $B$ to the real subspace $\cc/\cc\cap\hh_x\subset\gg/\hh_x$
is $\R$-valued and positive-definite. Indeed, if $\xi\in\cc$
then using the fact that $H$ is $\cc$-invariant we get
$$B(\xi,\xi)=H(\xi x,\xi x)\in\R.$$
Furthermore, this implies that $B(\xi,\xi)>0$ for $\xi\not\in\cc\cap\hh_x$.
Since the subspace $\cc/\cc\cap\hh$ generates $\gg/\hh_x$ over $\C$
this implies that the form $B$ on $\gg/\hh_x$ is non-degenerate.
In particular, the stabilizer subalgebra of $x^{\vee}$ is equal to $\hh_x$.
\ed

\begin{rem} The idea to look at vectors of minimal length on the
orbit goes back to the work of Kempf and Ness \cite{KN}.
One can give an alternative proof of proposition \ref{dual} using
the theorem of Mostow \cite{M} on self-adjoint groups.
\end{rem}

For nice representations the assertion of the above theorem holds also for
generic stabilizer subgroups.

\begin{prop}\label{dualnice} Assume that $\rho$ is
a nice representation. Then the dual representation $\rho^{\vee}$
is also nice. In this case
the generic stabilizer subgroups for $\rho$ and $\rho^{\vee}$ are
conjugate over $\ov{k}$.
\end{prop}

\Pf . The generic stabilizer subgroup for $\rho^{\vee}$ has form
$\sigma(\bH)$ where $\sigma$ is an automorphism of $\bG$
which restricts to the map $t\mapsto t^{-1}$ on
some maximal torus.
This immediately implies that $\rho^{\vee}$ is nice. The second
assertion follows from proposition \ref{dual}.
\ed

Let $Q$ denote the Killing form on $\gg$. The following lemma
is well-known but we include the proof since we couldn't
find the reference (in its statement $\hh$ can be replaced
by the Lie algebra of any reductive subgroup in $\bG$).

\begin{lem}\label{Killem}
The restriction of $Q$ to $\hh$ is non-degenerate.
\end{lem}

\Pf . We can assume that $k=\C$.
Let $C\subset\bG(\C)$ be a maximal compact subgroup containing
a maximal compact subgroup in $\bH(\C)$, and let $\cc$
be its Lie algebra. Then $\hh=\hh\cap\cc+i\hh\cap\cc$.
It remains to use the fact that the
restriction of $Q$ to $\cc$ is $\R$-valued and negative definite.
\ed

\begin{lem}\label{etale}
Let us denote $\bW=\bV^{\bH}$, $\bW_0=\bV_0\cap\bW$.

(a) For every point $x_0\in\bW_0$ there
exists a linear subspace $\bL\subset\bW_0$ and a Zariski open
neighborhood of zero $\bS\subset\bL$
such that the natural morphism
$$a:\bG/\bH\times (\bS+x_0)\ra\bV_0:(g\bH,x)\mapsto gx$$
is \'etale.\\

\noindent
(b) Let $\nu_{\bV}$ be a non-zero translation-invariant top-degree form on
$\bV$, $\om$ be a non-zero
$\bG$-invariant top-degree form on $\bG/\bH$. Then we have
$a^*\nu_{\bV}=\om\we\nu_{\bL}$ for some translation-invariant top-degree
form $\nu_{\bL}$ on $\bL$.
\end{lem}

\Pf . (a) Note that we have an isomorphism
$$\bG/\bH\times_{N(\bH)/\bH}\bW_0\ra\bV_0$$
where $N(\bH)$ is the normalizer of $\bH$ in $\bG$.
Let $\bL$ to be a complement in $\bW_0$
to the tangent space to $N(\bH)x_0$ at $x_0$. Then our assertion
follows from Luna's results in \cite{Lu} (in this situation $\bS+x_0$ is
an \'etale slice for the action of $N(\bH)/\bH$ on $\bW$). 

\noindent
(b) Since $\bG$ is semisimple and connected,
the form $\nu_{\bV}$ is $\bG$-invariant. Therefore, the pull-back
$a^*\nu_{\bV}$ is also $\bG$-invariant. On the other hand,
the map $a$ is linear in the second argument, hence
$a^*\nu_{\bV}$ is invariant with respect
to translations on $\bL$. This implies our assertion.
\ed


\begin{lem}\label{sscommlem}
Assume that $\bG$ is simple, $\rho$ is nice.
Then the generic stabilizer $\bH$ is either semisimple
or commutative.
\end{lem}

\Pf . This follows immediately from
Tables 1 and 2 of \cite{El1}. 
\ed

\subsection{Critical points}\label{critsec}

In this subsection we study critical points of the
restriction of a generic linear 
functional $x^{\vee}\in\bV^{\vee}_0$ to a generic orbit $\OO$
in a representation $\rho$ whose generic stabilizer subalgebra
is reductive. Recall that a critical point $x$ of a function
$\phi$ is called {\it non-degenerate} if the quadratic form of second
derivatives of $\phi$ at $x$ is non-degenerate. 

\begin{lem}\label{critlem}
Let $x$ be a critical point of $x^{\vee}|_{\OO}$. Then:

\noindent (a)
$x$ has the same stabilizer subalgebra as $x^{\vee}$;\\

\noindent (b)
$x$ is non-degenerate if and only if $\gg x\cap (\gg x^{\vee})^{\perp}=0$.
\end{lem}

\Pf . (a) Without loss of generality we can assume that
$x^{\vee}\in\bW^{\vee}_0$.
The condition that $x$ is a critical point of
$x^{\vee}|_{\OO}$ is equivalent to $x^{\vee}(\gg x)=0$, i.e.,
$x\in (\gg x^{\vee})^{\perp}$.
On the other hand, by lemma \ref{sumlem} the group
$\bH$ acts trivially on $\bV^{\vee}/\gg x^{\vee}$. It follows that
$(\gg x^{\vee})^{\perp}$ is contained in $\bV^{\bH}=\bW$, hence,
$x$ is stabilized by $\bH$.

\noindent
(b) The second derivative of $x^{\vee}|_{\OO}$
at a critical point $x$ is the following
symmetric bilinear form on the tangent space $T_x\OO=\gg x$:
$B_{x,x^{\vee}}(\xi_1 x,\xi_2 x)=\lan x^{\vee}, \xi_1\xi_2 x\ran$,
where $\xi_1,\xi_2\in\gg$.
The kernel of this form is $\gg x\cap (\gg x^{\vee})^{\perp}$.
\ed

Let $\bI\subset\bV_0\times\bV_0^{\vee}$
be the subvariety consisting of $(x,x^{\vee})$ such that
$x$ is a critical point of $x^{\vee}|_{\OO_x}$. 
Let us consider the natural morphism
$$\bff:\bI\ra\bV_0/\bG\times\bV_0^{\vee}.$$
By the definition, $\bff^{-1}(\OO,x^{\vee})$ consists
of critical points of $x^{\vee}|_{\OO}$.

\begin{prop}\label{critprop}
The morphism $\bff$ is dominant
and there exists a non-empty $\bG$-invariant
open subset $\bU\subset \bV_0/\bG\times
\bV_0^{\vee}$ (where $\bG$ acts on the second factor) such that $\bff$
is \'etale over $\bU$. For every pair $(\OO,x^{\vee})$
in $\bU$ 
the set of critical points of $x^{\vee}|_\OO$ is finite and non-empty.
Furthermore, all these critical points
are non-degenerate.
\end{prop}

\Pf . By the definition, the subvariety $\bI\subset\bV_0\times\bV_0^{\vee}$
consists of $(x,x^{\vee})$ such that
$\lan x^{\vee},\gg x\ran=0$. Let $p_1:\bI\ra\bV_0$, $p_2:\bI\ra\bV_0^{\vee}$
be the natural projections.
Note that the fibers of both these projections
are open subsets in linear spaces of dimension
$\dim \bV-\dim\OO=\dim\bV/\bG$. Hence,
$\bI$ is smooth and irreducible
of dimension $\dim\bV/\bG+\dim\bV$ (if non-empty).
Thus, it suffices to prove that the morphism $\bff$
is dominant. We can assume that $k=\C$. Let $(x,x^{\vee})$ be
a pair constructed in the proof of proposition \ref{dual}.
Then $(x,x^{\vee})\in\bI$. 
We claim that the tangent map to $\bff$ at $(x,x^{\vee})$ is an
isomorphism. Indeed, the relative tangent space of the projection
$p_2:\bI\ra\bV_0^{\vee}$ at $(x,x^{\vee})$ is $(\gg x^{\vee})^{\perp}$.  
Therefore, we have to check that the linear map
$(\gg x^{\vee})^{\perp}\ra \bV/\gg x$ is an isomorphism.
It remains to note that the non-degeneracy
of the form $B$ constructed in the proof of proposition \ref{dual}
implies that $\gg x\cap (\gg x^{\vee})^{\perp}=0$.
This proves our claim. The last assertion of the proposition follows
from lemma \ref{critlem} (b) and from the fact that $\bff$ is generically
\'etale.
\ed

\subsection{Spin-coverings}\label{spinsec}

In this subsection we are going to define certain double coverings of generic
stabilizers which will play an important role later. We assume
that $\rho$ is a nice representation.

Let $x\in\bV_0$, $x^{\vee}\in\bV^{\vee}_0$ be
points such that $x$ is a non-degenerate critical point
of $x^{\vee}|_{\OO_x}$. Let $B_{x,x^{\vee}}$ be the corresponding
quadratic form of second derivatives
of $x^{\vee}|_{\OO_x}$ on the tangent space $T_x:=T_x\OO_x$. Then
the action of the stabilizer $\bH_x$ of $x$ on $T_x$ preserves
$B_{x,x^{\vee}}$, hence, we get a homomorphism  
$\bH_x\ra \Aut(T_x,B_{x,x^{\vee}}).$
Since $\bH_x$ is connected, we obtain a homomorphism to the
corresponding special orthogonal group
$$\iota_{x,x^{\vee}}:\bH_x\ra\SO(T_x,B_{x,x^{\vee}}).$$

Recall that for every vector space $T$ equipped with a non-degenerate
quadratic form $B$
one can define the spin-covering $\Spin(T,B)\ra\SO(T,B)$. 
It is defined in the standard way when
$\dim T\ge 3$ (see e.g. \cite{Bou}).
If $\dim T<3$ then we define the spin-covering
as the restriction of the standard spin-covering of $\SO(T\oplus U)$,
where $U$ is an orthogonal space of large dimension.
The group $\Spin(T,B)$ is a central extension of $\SO(T,B)$
by $\{\pm 1\}$. In the case $\dim(T)<3$ this extension is trivial unless
$\dim T=2$ and the quadratic form $B$ is anisotropic. In the latter
case $\SO(T,B)$ is a $1$-dimensional torus of the form
$R^{(1)}_{k'/k}(\G_m)$ for some quadratic extension $k\subset k'$,
$\Spin(T,B)$ is isomorphic to $\SO(T,B)$ and the map
$\Spin(T,B)\ra\SO(T,B)$ is identified with
the map $t\mapsto t^2$.

In proposition \ref{spinprop} below we compute 
the pull-back of the spin-covering
$\Spin(T_x,B_{x,x^{\vee}})\ra\SO(T_x,B_{x,x^{\vee}})$
by $\iota_{x,x^{\vee}}$.
Recall that by lemma \ref{Killem}
the restriction of $Q$ to $\hh_x$ is non-degenerate.
Consider the homomorphism
$$\iota_x:\bH_x\ra\SO(\hh_x,Q|_{\hh_x})$$
induced by the adjoint action of $\bH_x$.
Let $\wt{\bH}_x\ra\bH_x$ be the pull-back of the spin-covering by
$\iota_x$. Note that if $\bH_x$ is commutative then $\iota_x$
is trivial, so the extension $\wt{\bH}_x\ra\bH_x$ splits in this case.

\begin{prop}\label{spinprop}
There is a unique isomorphism of the following two
central extensions of $\bH_x$ by $\{\pm 1\}$:
the pull-back of the spin-covering by $\iota_{x,x^{\vee}}$
and $\wt{\bH}_x$.
\end{prop}

\Pf . We can reformulate the statement as follows: there exists
a unique lifting of $\iota_{x,x^{\vee}}$ to a
homomorphism $\wt{\bH}_x\ra\Spin(T_x,B_{x,x^{\vee}})$
which maps
$\{\pm 1\}\subset\wt{\bH}_x$ to $\{\pm 1\}\subset\Spin(T_x,B_{x,x^{\vee}})$
identically.

First, let us prove the uniqueness. Indeed, two such liftings
differ by a homomorphism $\bH_x\ra\{\pm 1\}$. Since $\bH_x$ is connected
such a homomorphism should be trivial. Now by uniqueness it suffices
to prove the existence of a lifting over $\ov{k}$. Therefore,
we can assume that $k=\C$.

Assume first that $\dim(T_x)\ge 3$ and
the twofold covering $\wt{\bH}_x\ra\bH_x$ is non-trivial. Then $\wt{\bH}_x$
is connected. We claim that the homomorphism
\begin{equation}\label{homspin}
\pi_1(\bH_x)\ra\pi_1(\SO(T_x,B_{x,x^{\vee}}))=\{\pm 1\}
\end{equation}
induced by $\iota_{x,x^{\vee}}$ is surjective with the kernel
$\pi_1(\wt{\bH}_x)\subset\pi_1(\bH_x)$.
Indeed, let $X$ be the space of
$\bH_x$-invariant non-degenerate symmetric forms on $T_x$.
We have a continuous family of homomorphisms $\bH_x\ra\SO(T_x,B)$
parametrized by $B\in X$. Since $X$ is connected,
the induced homomorphism of fundamental groups
$\pi_1(\bH_x)\ra\pi_1(\SO(T_x,B))$ 
does not depend on $B$.
Therefore, in the above claim we can replace $B_{x,x^{\vee}}$
by any other form in $X$.
Now using the natural identification of
$T_x$ with $\gg/\hh_x=\hh_x^{\perp}\subset \gg$ (the orthogonal
complement to $\hh_x$ in $\gg$ with respect to the Killing form $Q$)
we take $B=Q|_{\hh_x^{\perp}}$. Consider the following commutative
diagram of homomorphisms
\begin{equation}
\begin{array}{ccc}
\bH_x &\lrar{} &\SO(\hh_x,Q|_{\hh_x})\times
\SO(\hh_x^{\perp},Q|_{\hh_x^{\perp}})\\
\ldar{} & & \ldar{}\\
\bG &\lrar{} &\SO(\gg,Q)\
\end{array}
\end{equation}
This diagram implies that the following composition
$$\pi_1(\bH_x)\ra\pi_1(\SO(\hh_x,Q|_{\hh_x}))\times
\pi_1(\SO(\hh_x^{\perp},Q|_{\hh_x^{\perp}}))\ra\pi_1(\SO(\gg,Q))$$
is trivial (since it factors through $\pi_1(\bG)=1$).
The assumption that the covering $\wt{\bH}_x\ra\bH_x$
is non-trivial implies that $\bH_x$ has a non-trivial semisimple component
(in particular, $\dim\hh_x\ge 3$) and that the map
$\pi_1(\bH_x)\ra\pi_1(\SO(\hh_x,Q|_{\hh_x}))$ is surjective with
the kernel $\pi_1(\wt{\bH}_x)\subset\pi_1(\bH_x)$.
On the other hand, both the maps
$$\pi_1(\SO(\hh_x,Q|_{\hh_x}))\ra\pi_1(\SO(\gg,Q))$$
and
$$\pi_1(\SO(\hh_x^{\perp},Q|_{\hh_x^{\perp}}))\ra\pi_1(\SO(\gg,Q))$$
are isomorphisms.
Our claim immediately follows from this.
This implies that the homomorphism
$\wt{\bH}_x\ra\SO(T_x,B_{x,x^{\vee}})$ lifts to a homomorphism
$\wt{\bH}_x\ra\Spin(T_x,B_{x,x^{\vee}})$. It is easy to see that the
latter homomorphism is non-trivial on $\{\pm 1\}\subset\wt{\bH}_x$.
Indeed, otherwise it would factor through $\bH_x$ which contradicts
to non-triviality of the homomorphism (\ref{homspin}). This finishes
the proof in this case.

Now assume that the covering $\wt{\bH}_x\ra\bH_x$ is trivial (and
$\dim(T_x)\ge 3$).
Then the map $\pi_1(\bH_x)\ra\pi_1(\SO(\hh_x,Q|_{\hh_x}))$ is
trivial (when $\dim\hh_x\le 2$ this follows from commutativity of
$\bH_x$). As above we deduce from this that
the map $\pi_1(\bH_x)\ra\pi_1(T_x,B_{x,x^{\vee}})$
is also trivial. Therefore, the homomorphism
$\bH_x\ra\SO(T_x,B_{x,x^{\vee}})$ factors through
$\Spin(T_x,B_{x,x^{\vee}})$.

Finally, if $\dim(T_x)<3$ then we replace the orthogonal
space $(T_x,B_{x,x^{\vee}})$ by its direct sum with a fixed orthogonal
space of large dimension (on which $\bG$ acts trivially) and 
apply the same argument as above. 
\ed

In the case when $k=E$ is a local field, the exact sequence
$1\ra\{\pm 1\}\ra \wt{\bH}_x\ra\bH_x\ra 1$ 
gives a map of Galois cohomologies
\begin{equation}\label{signfnc}
\eps_{\bH_x}:H^1(E,\bH_x)\ra H^2(E,\{\pm 1\})\simeq\{\pm 1\}.
\end{equation}
This is the sign function which will play an important role below.

\section{Nice representations over local fields}\label{localsec}

In this section $\rho:\bG\ra\Aut(\bV)$ denotes
a nice representation over a local field $E$.
We formulate our main theorems about the behaviour of
$\bG$-stable functions  and distributions
on $\bV(E)$ under the Fourier transform.
Also, we analyze the Fourier transform of certain stable
functions using the stationary phase approximation 
(in the case when $E$ is $p$-adic).

\subsection{$\bG$-inner forms and stable $\bG$-equivalence}\label{inner}

Let $x$ be a point in $\bV_0(E)$.
Note that since $\bH_x$ is reductive, 
the set $H^1(E,\bH_x)$ is finite (see e.g. \cite{PR}).
Consider the subset $P_x\subset\bG(\ov{E})$ consisting of the
elements $g$ such that $g^{-1}\sigma(g)\in\bH_x$
for all $\sigma\in\Ga_E$. For every $g\in P_x$ the $1$-cocycle
$$e_g(\sigma)=g^{-1}\sigma(g)$$
gives a class in $\bH_x$. It is easy to see that $P_x$ is a union of 
right $\bH_x(\ov{E})$-cosets and of left $\bG(E)$-cosets, and
the assignment $g\mapsto e_g$ defines a bijection
$$\bG(E)\backslash P_x/\bH_x(\ov{E})\wt{\ra} 
\ker(H^1(E,\bH_x)\ra H^1(E,\bG)).$$
On the other hand, we have a natural bijection  
$$P_x/\bH_x(\ov{E})\wt{\ra}\OO_x(E): g\mapsto gx.$$
In particular, we can identify the set of $\bG(E)$-orbits on 
$\OO_x(E)$ with $\ker(H^1(E,\bH_x)\ra H^1(E,\bG))$.
In the case when $E$ is $p$-adic, we have $H^1(E,\bG)=0$ 
since $\bG$ is simply connected (see \cite{PR}). Therefore, in this
case we have a bijection between $\bG(E)\backslash\OO_x(E)$ and
$H^1(E,\bH_x)$.

Let us call two points $x,x'\in\bV(E)$ {\it stably $\bG$-equivalent}
if there exists an element $g\in\bG(\ov{E})$ such that $gx=x'$, i.e.,
if $\OO_x=\OO_{x'}$.
Let $\bK,\bK'\subset \bG$ be subgroups defined over $E$. Let us say
that $\bK'$ is a $\bG$-{\it inner form} of $\bK$ if there exists
an element $g\in\bG(\ov{E})$ such that $g\bK(\ov{E}) g^{-1}=\bK'(\ov{E})$ and
$g^{-1}\sigma(g)\in\bK(\ov{E})$ for every $\sigma\in\Ga_E$.
It is easy to see that this defines an equivalence relation between
subgroups of $\bG$ defined over $E$. This definition is motivated by
the following lemma.

\begin{lem}\label{Ginner}
Let $x\in\bV(E)$.\newline
(a) If $x'\in\bV(E)$ is stably $\bG$-equivalent to $x$ then
$\bH_{x'}$ is a $\bG$-inner form of $\bH_x$.\newline
(b) Conversely, if $\bH'$ is a $\bG$-inner form of $\bH_x$ 
then $\bH'$ is the stabilizer of
some $E$-point which is stably $\bG$-equivalent to $x$.\newline
(c) The image of the natural map
$$a:\bG/\bH(E)\times \bW_0(E)\ra\bV_0(E)$$
consists of all points in $\bV_0(E)$ whose stabilizer is a $\bG$-inner
form of $\bH$.
\end{lem}

The proof is straightforward.

\begin{rem} It is easy to see that if $\bK'$ is a $\bG$-inner form of 
$\bK$ then $\bK'$ is an inner form of $\bK$ in the usual sense. 
On the other hand, if $\bK$ is semisimple and simply connected,
$E$ is $p$-adic, then $\bK'$ is a $\bG$-inner form of $\bK$ if and only
if $\bK'$ is conjugate to $\bK$ over $E$.
\end{rem}

Let $x$ and $x'$ be a pair of stably $\bG$-equivalent points in
$\bV_0(E)$ and let $\OO=\OO_x=\OO_{x'}$ be the corresponding orbit.
Then $\bG(E)$-orbits on $\OO(E)$ can be identified with a subset in
$H^1(E,\bH_x)$ and with a subset in $H^1(E,\bH_{x'})$. On the other hand,
$\bH_{x'}$ is obtained from $\bH_x$ by twisting with a cohomology class
in $H^1(E,\bH_x)$, so we have a canonical
identification $H^1(E,\bH_{x'})\simeq H^1(E,\bH_x)$. It is easy
to see that these three identifications are compatible.


\subsection{Local Kottwitz invariant and local sign function}
\label{localKotsec}

For every connected reductive group $\bH$ over a field $k$ let us denote by
$Z(\hat{\bH})$ the centre of the Langlands dual group
(equipped with an action of the Galois group $\Ga_k$). Note that
$\hat{\bH}$ is defined canonically up to an inner conjugation, hence,
$Z(\hat{\bH})$ is defined canonically and carries an action of
$\Ga_k$.
Furthermore, an isomorphism $i:\bH\ra\bH'$ over $\ov{k}$,
such that $i^{-1}\sigma(i)$
is inner for all $\sigma\in\Ga_k$ (an {\it inner twisting}),
induces an isomorphism
of $Z(\hat{\bH})$ with $Z(\hat{\bH'})$ as $\Ga_k$-modules.
Following Kottwitz we define
$$A(\bH/k):=\pi_0(Z(\hat{\bH})^{\Ga_k})^D$$
where for a finite group $A$ we denote by $A^D$ the dual group.
When $k=E$ is a local field, the local Kottwitz invariant is a functorial map
$$\inv=\inv_E:H^1(E,\bH)\ra A(\bH/E)$$
constructed in \cite{K2}.
The definition of $\inv$ is a generalization
of the isomorphism derived from Tate-Nakayama duality in the
case when $\bH=\bT$ is a torus. Indeed, this duality gives an
isomorphism $H^1(E,\bT)\simeq H^1(E,X^*(\bT))^D$, where
$X^*(\bT)$ is the module of characters of $\bT$. Now from
the exact sequence of Galois modules
$$0\ra X^*(\bT)\ra X^*(\bT)\otimes\C\ra X^*(\bT)\otimes\C^*\ra 0$$
one gets an isomorphism
$$H^1(E,X^*(\bT))\simeq\coker((X^*(\bT)\otimes\C)^{\Ga_E}\ra
(X^*(\bT)\otimes\C^*)^{\Ga_E}).$$
The latter group can be immediately identified with
$\pi_0((X^*(\bT)\otimes\C^*)^{\Ga_E})$.
In the general case (when $\bH$ is not necessarily a torus),
Kottwitz showed in \cite{K2} that
$\inv_E$ is an isomorphism for $p$-adic $E$.
In particular, for such $E$ we obtain the structure of abelian
group on $H^1(E,\bH)$.

For any $1$-cocycle $e:\Ga_E\ra\bH(\ov{E})$
we can consider the $E$-group $\bH^e$ obtained from $\bH$ by
inner twisting with $e$. By definition $\bH^e(\ov{E})=\bH(\ov{E})$
while the action of an element $\sigma\in\Ga_E$
on $\bH^e(\ov{E})$ differs from its action on $\bH(\ov{E})$ by the
inner automorphism associated with $e(\sigma)$. In particular,
we have $Z(\hat{\bH^e})=Z(\hat{\bH})$, hence $A(\bH^e/E)=A(\bH/E)$.
According to lemma 1.4 of \cite{K2},
the following diagram is commutative
\begin{equation}
\begin{array}{ccc}
H^1(E,\bH^e) &\lrar{\inv} & A(\bH/E)\\
\ldar{i_e}& & \ldar{t_{\inv(e)}}\\
H^1(E,\bH)&\lrar{\inv} & A(\bH^e/E)=A(\bH/E)
\end{array}
\end{equation}
where $t_{\inv(e)}:A(\bH/E)\ra A(\bH/E)$ is the translation
by $\inv(e)\in A(\bH/E)$, $i_e:H^1(E,\bH^e)\ra H^1(E,\bH)$
is the canonical identification induced by $e$.
Thus, in the case when $E$ is $p$-adic,
the isomorphism $i_e$ does not respect
the group structures on the sets $H^1(E,\bH^e)$ and $H^1(E,\bH)$,
but rather respects the structures of
principal homogeneous spaces over $A(\bH/E)=A(\bH^e/E)$.

Now assume that $\OO\subset\bV_0$ is an orbit. Then
we have a system of compatible isomorphisms between
the groups $A(\bH_x/E)$ for $x\in \OO(E)$.
\footnote{For two points $x,x'\in\OO(E)$
such that $\bH_x=\bH_{x'}$, the corresponding isomorphism
$A(\bH_x/E)\ra A(\bH_{x'}/E)$ is not necessarily the identity:
it corresponds to the action of some element in the normalizer
of $\bH_x$.}
Let us denote the corresponding group isomorphic to
all $A(\bH_x/E)$ by $A(\OO/E)$.
Assume for a moment that $E$ is $p$-adic.
Then the set of $\bG(E)$-orbits
on $\OO(E)$ has a natural structure of a principal homogeneous
space over $A(\OO/E)$. Thus, for every pair of points $x,x'\in\OO(E)$
we can define an element $\inv(x,x')\in A(\OO/E)$ such that
$\bG(E)x'$ is obtained from $\bG(E)x$ by the action of $\inv(x,x')$.
This definition extends to the case of archimedian $E$ as follows:
$$\inv(x,gx)=\inv(e_g)$$
where $g\in P_x$, $e_g$ is the corresponding cohomology class in
$H^1(E,\bH_x)$, $\inv(e_g)$ is its local Kottwitz invariant in
$A(\bH_x/E)\simeq A(\OO/E)$.
It is easy to see that the following properties are satisfied:
\begin{equation}\label{invprop}
\begin{array}{l}
\inv(x,x')+\inv(x',x)=0,\\
\inv(x,x')+\inv(x',x'')=\inv(x,x''),\\
\inv(tx,tx')=\inv(x,x')
\end{array}
\end{equation}
where $g\in\bG(\ov{E})$, $t\in E^*$,
$e_g$ is the $1$-cocycle of $\Ga_E$ with values in $\bH_x$
defined above.

Let us fix a point $x\in\bV_0(E)$ and set $\bH=\bH_x$.
Recall that in section \ref{spinsec} we have defined a
map
$$\eps_{\bH}:H^1(E,\bH)\ra H^2(E,\{\pm 1\})=\{\pm 1\}$$
induced by the central extension $1\ra\{\pm 1\}\ra\wt{\bH}\ra\bH\ra 1$ 
(the pull-back of the spin covering associated with $Q|_{\hh_x}$).
We are going to construct a character
$$\sign=\sign_{\bH}:A(\bH/E)\ra\{\pm 1\}$$ 
such that $\eps_{\bH}=\sign_{\bH}\circ\inv_E$. 
For this we note that the above central extension is induced
by the similar extension $1\ra\{\pm 1\}\ra\wt{\bH}_{ad}\ra\bH_{ad}\ra 1$ of
the adjoint group $\bH_{ad}=\bH/Z(\bH)$. Therefore, by
functoriality of $\inv_E$, it suffices to construct the
character $\sign_{\bH}$ in the case when $\bH$ is adjoint.
Thus, we can assume that $\bH$ is semisimple.
Let $u:\bH_{sc}\ra\bH$ be the universal covering of $\bH$ (it is defined
over $E$), $\bC$ be the kernel of $u$. Then the homomorphism $u$ lifts 
uniquely to a homomorphism $\wt{u}:\bH_{sc}\ra\wt{\bH}$.
Now the restriction of $\wt{u}$ to $\bC$ gives a homomorphism
$\chi:\bC\ra\{\pm 1\}$ defined over $E$. We can consider $\chi$ as
an element of order $2$ in $X^*(\bC)^{\Ga_E}$.
It remains to notice that there is an isomorphism
$X^*(\bC)\simeq Z(\hat{\bH})$ of
$\Ga_E$-modules, so we can consider $\chi$ as
a character $\sign=\sign_{\bH}:A(\bH/E)\ra\{\pm 1\}$. 
The following result shows that this is the character we were looking for.

\begin{lem}\label{signinv} One has the following equality of maps
from $H^1(E,\bH)$ to $\{\pm 1\}$:
$$\eps_{\bH}=\sign_{\bH}\circ\inv_E.$$
\end{lem}

\Pf . By duality for finite groups we have an isomorphism
$$H^2(E,\bC)\simeq H^0(E,X^*(\bC))^D\simeq A(\bH/E).$$
According to Lemma 1.8 of \cite{K2},
under this isomorphism the map
$\inv_E:H^1(E,\bH)\ra A(\bH/E)$ 
can be identified with the map $H^1(E,\bH)\ra H^2(E,\bC)$ coming
from the exact sequence $1\ra\bC\ra\bH_{sc}\ra\bH\ra 1$.
On the other hand, the character of $H^2(E,\bC)$ corresponding to $\sign_{\bH}$
is the homomorphism on $H^2$ induced by the homomorphism $\chi:\bC\ra\{\pm 1\}$.
Therefore, the composition $\sign_{\bH}\circ\inv_E$ coincides with
the map $H^1(E,\bH)\ra H^2(E,\{\pm 1\})$ coming from the exact
sequence $1\ra\{\pm 1\} \ra\wt{\bH}\ra\bH\ra 1$, which is the definition of
$\eps_{\bH}$.
\ed

Recall that for every pair of points $x,x'\in\OO(E)$, where $\OO\subset\bV_0$
is a $\bG$-orbit, we have defined the sign 
$\eps(x,x')=\eps(Q|_{\hh_x},Q|_{\hh_x'})$.
It is well-known that if a quadratic form $B'$ is obtained from
a non-degenerate quadratic form $B$ by the twist with an element
$e\in H^1(E,\SO(B))$ then the relative Hasse-Witt invariant
$\eps(B,B')$ is equal to the image of $e$ under the coboundary homomorphism
$H^1(E,\SO(B))\ra H^2(E,\{\pm 1\})=\{\pm 1\}$ coming from the spin-covering
(see e.g. \cite{Spr1}). 
This implies the following relation between $\eps$ and
the sign function (\ref{signfnc}) defined in \ref{spinsec}:
\begin{equation}\label{signrelation}
\eps(x,gx)=\eps_{\bH_x}(e_g)
\end{equation}
where $x\in\bV_0(E)$, $g\in P_x$.
Comparing the definition of $\inv(\cdot,\cdot)$ with (\ref{signrelation})
and using lemma \ref{signinv} we get the following formula:
$$\eps(x,x')=\sign_{\bH_x}(\inv(x,x')).$$

\subsection{Critical points and stable equivalence}

Assume that we have a pair of points $x\in\bV_0(E)$,
$x^{\vee}\in\bV_0^{\vee}(E)$ such that $x$ is a critical point
of $x^{\vee}|_{\OO_x}$.
Proposition \ref{spinprop} implies that
the following diagram is commutative
\begin{equation}\label{comdiaghasse}
\begin{array}{ccc}
H^1(E,\bH_x) &\lrar{\eps_{\bH_x}}& \{\pm 1\}\\
\ldar{H^1(\iota_{x,x^{\vee}})} & & \ldar{\id}\\
H^1(E,\SO(T_x,B_{x,x^{\vee}})) &\lrar{}&\{\pm 1\}
\end{array}
\end{equation}
where the lower horizontal arrow is the coboundary homomorphism
associated with the spin-covering of $\SO(T_x,B_{x,x^{\vee}})$.

This observation leads to the following result.

\begin{lem}\label{quadr}
Let $\om$ be a $\bG$-invariant top-degree form on $\OO_x$.
Recall that for every $g\in P_x$ one has $gx\in\bV_0(E)$ and
$gx^{\vee}\in\bV_0^{\vee}(E)$.

\noindent
(a) One has $\det(B_{x,x^{\vee}})/\om_x^2=
\det(B_{gx,gx^{\vee}})/\om_{gx}^2$.

\noindent
(b) The quadratic forms $B_{x,x^{\vee}}$ and $B_{gx,gx^{\vee}}$
have the same determinant modulo squares. Their 
relative Hasse-Witt invariant is given by
$$\eps(B_{gx,gx^{\vee}},B_{x,x^{\vee}})=\eps(x,gx)=\eps_{\bH_x}(e_g).$$
\end{lem}

\Pf . (a) This follows from the $\ov{E}$-isomorphism of data
$(T_v,B_{v,v^{\vee}},\om_v)$ and $(T_{gv}, B_{gv,gv^{\vee}},\om_{gv})$
given by the action of $g$.

\noindent
(b) Let $e_g\in H^1(E,\bH_v)$ be the cohomology class defined by
the cocycle $\sigma\mapsto g^{-1}\sigma(g)$. Then the quadratic form
$B_{gv,gv^{\vee}}$ is obtained from $B_{v,v^{\vee}}$ by twisting
with the class $H^1(\iota_{v,v^{\vee}})(e_g)\in
H^1(E,\SO(T_v,B_{v,v^{\vee}}))$. Now the diagram (\ref{comdiaghasse}) 
implies that
$\eps(B_{gx,gx^{\vee}},B_{x,x^{\vee}})=\eps_{\bH_x}(e_g)$.
\ed

\subsection{Stable and antistable functions and distributions}
\label{stabledistrsec}

Let us denote by $\SS(\bV(E))_{\bG(E)}$ the space of $\bG(E)$-coinvariants
in $\SS(\bV(E))$. We have the natural projection
$\SS(\bV(E))\ra\SS(\bV(E))_{\bG(E)}:\phi\mapsto\ov{\phi}$.
On the other hand, for every
$\phi\in\SS(\bV(E))$ we can define a function $I(\phi)$ on
$\bG(E)\backslash \bV_0(E)$
by the formula
$$I(\phi)(y)=\int_{x\in\bG(E)y}\phi(x)|\om_y|$$
where $\om_y$ is a $\bG$-invariant top-degree form on $\OO_y$ (the integral
is convergent since the orbit $\OO_y$ is closed in $\bV$).
It is clear that $I(\phi)$ depends only on $\ov{\phi}$, so we will denote
$I(\ov{\phi})=I(\phi)$. Although we will not need this fact,
it is worth mentioning that for a pair of
functions $\phi,\phi'\in\SS(\bV_0(E))$ one has $\ov{\phi}=\ov{\phi'}$
if and only if
$I(\phi)=I(\phi')$ (see \cite{GK}).

\noindent {\bf Definition}.

\noindent (i)
An element $\ov{\phi}\in\SS(\bV(E))_{\bG(E)}$ is called {\it stable}
if for every $\bG$-orbit $\OO\subset\bV_0$, the restriction
of the function $I(\ov{\phi})$ to $\bG(E)\backslash\OO(E)$ is constant. 
We denote by $\SS(\bV(E))^{st}\subset\SS(\bV(E))_{\bG(E)}$
the subspace of stable elements.
Similarly, we define a subspace of {\it antistable} elements
$\SS(\bV(E))^{as}\subset\SS(\bV(E))_{\bG(E)}$. By the definition,
an element $\ov{\phi}\in\SS(\bV(E))_{\bG(E)}$ is antistable
if for every $\bG$-orbit $\OO\subset\bV_0$ the total sum
of the function $I(\ov{\phi})$ over $\bG(E)\backslash\OO(E)$
is zero.

\noindent (ii)
Let $\DD(\bV(E))$ denote the space of distributions on $\bV(E)$,
i.e., functionals on $\SS(\bV(E))$.
Note that a
$\bG(E)$-invariant distribution $\a\in\DD(\bV(E))$ descends to
a functional on $\SS(\bV(E))_{\bG(E)}$. Now a $\bG(E)$-invariant
distribution $\a$ is called {\it stable}
(resp. {\it antistable}) if $\a(\SS(\bV(E))^{as})=0$
(resp. $\a(\SS(\bV(E))^{st})=0$).
We denote by $\DD(\bV(E))^{st}\subset\DD(\bV(E))$
(resp. $\DD(\bV(E))^{as}\subset\DD(\bV(E))$) the subspace
of stable (resp. antistable) distributions.

\vspace{2mm}

For a $\bG$-orbit $\OO\subset\bV_0$ 
and a non-zero $\bG$-invariant top-degree form $\om$ on $\OO$
defined over $E$, we can define a stable distribution
$\delta_{\OO,\om}\in\DD(\bV(E))$ by the formula
$$\delta_{\OO,\om}(\phi)=\int_{\OO(E)}\phi |\om|$$
where $\phi\in\SS(\bV(E))$.
When $\OO(E)=\emptyset$ we set $\delta_{\OO,\om}=0$.
When the choice of $\om$ is clear or is not important we will
abbreviate $\delta_{\OO,\om}$ to $\delta_{\OO}$. 
By the definition,
an element $\ov{\phi}\in\SS(\bV(E))_{\bG(E)}$
is antistable if and only if
$\delta_{\OO}(\ov{\phi})=0$
for all $\bG$-orbits $\OO\subset\bV_0$.

More generally, for every point $x\in\bV_0(E)$ and for a character
$\kappa:A(\OO_x/E)\ra\C^*$ we can define a $\bG(E)$-invariant
distribution
$\delta^{\kappa}_{\OO_x}=\delta^{\kappa}_{\OO_x,\om}$ by the formula
$$\delta^{\kappa}_{\OO_x,\om}(\phi)=
\int_{x'\in\OO_x(E)}\kappa(\inv(x,x'))\phi(x')
|\om(x')|.$$
Note that if $y$ is stably equivalent to $x$ then
$$\delta^{\kappa}_{\OO_y,\om}=
\kappa(\inv(y,x))\cdot\delta^{\kappa}_{\OO_x,\om}.$$
In the case when $E$ is $p$-adic,
an element $\ov{\phi}\in\SS(\bV(E))_{\bG(E)}$ is stable
if and only if $\delta^{\kappa}_{\OO_x}(\ov{\phi})=0$
for all $x\in\bV_0(E)$ and all non-trivial characters $\kappa$ of
$A(\OO_x/E)$.

Also, for every $E$-orbit $\OO\subset\bV_0$ we define the distribution
$$\delta^{\eps}_{\OO,\om}(\phi)=\int_{y\in\OO(E)}\eps(y)\phi(y)|\om(y)|$$
where $\eps:\bV_0(E)\ra\{\pm 1\}$ is the function defined in the introduction.
The results of section \ref{localKotsec}
show that this distribution corresponds to some character of $A(\OO/E)$
as above. More precisely, for any point $x\in\bV_0(E)$ we have
$$\delta^{\eps}_{\OO_x,\om}=\eps(x)\cdot\delta^{\sign_{\bH_x}}_{\OO_x,\om}.$$

The Fourier transform (associated with some choice of a non-trivial
additive character $\psi$) induces a well-defined operator
$$\FF:\SS(\bV(E))_{\bG(E)}\ra\SS(\bV^{\vee}(E))_{\bG(E)}.$$
We want to describe the images of the subspaces
$\SS(\bV(E))^{st}$ and $\SS(\bV(E))^{as}$ under $\FF$.
This is equivalent to describing the images of the spaces of distributions
$\DD(\bV(E))^{st}$ and $\DD(\bV(E))^{as}$.

\noindent {\bf Definition}.

\noindent (i)
An element $\ov{\phi}\in\SS(\bV(E))_{\bG(E)}$ is called $\eps$-{\it stable}
if for every orbit $\OO\subset\bV_0$ and for every pair of 
points $y,y'\in\OO(E)$, one has
$I(\ov{\phi})(y')=\eps(y,y')I(\ov{\phi})(y)$.
We denote by $\SS(\bV(E))^{st}_{\eps}$ 
the subspace of $\eps$-stable elements in $\SS(\bV(E))_{\bG(E)}$.
Similarly we define the subspace
$\SS(\bV(E))^{as}_{\eps}\subset\SS(\bV(E))_{\bG(E)}$
of $\eps$-{\it antistable} elements. By the definition,
an element $\ov{\phi}\in\SS(\bV(E))_{\bG(E)}$ is $\eps$-antistable 
if for every $y\in\bV_0(E)$ one has
$\sum_{y'\in\bG(E)\backslash\OO_y(E)}\eps(y,y')I(\ov{\phi})(y')=0$.

\noindent (ii)
Dually, we define the
subspace $\DD(\bV(E))^{st}_{\eps}\subset\DD(\bV(E))^{\bG(E)}$
(resp. $\DD(\bV(E))^{as}_{\eps}\subset\DD(\bV(E))^{\bG(E)}$) 
of $\eps$-stable (resp. $\eps$-antistable) distributions, so that
$\DD(\bV(E))^{st}_{\eps}$ is the annihilator of $\SS(\bV(E))^{as}_{\eps}$
(resp. $\DD(\bV(E))^{as}_{\eps}$ is the annihilator of
$\SS(\bV(E))^{st}_{\eps}$).
\vspace{2mm}

Note that the distributions $\delta^{\eps}_{\OO}$ are $\eps$-stable,
and an element $\ov{\phi}\in\SS(\bV(E))$ is $\eps$-antistable if and only
if $\delta^{\eps}_{\OO}(\ov{\phi})=0$ for all $\bG$-orbits
$\OO\subset\bV_0$.

We say that a function $\phi\in\SS(\bV(E))$ is stable (resp. antistable,
$\eps$-stable, $\eps$-antistable) if this is true for the corresponding
element $\ov{\phi}\in\SS(\bV(E))_{\bG(E)}$.
Clearly, a function $\phi\in\SS(\bV_0(E))$
is $\eps$-stable (resp. $\eps$-antistable)
if and only if $\eps\cdot \phi$ is stable (resp. antistable).

Our main result is the following theorem.

\begin{thm}\label{mainthm}
Let $\rho:\bG\ra\GL(\bV)$ be a nice representation
of a simply connected semisimple group $\bG$ over a $p$-adic field
$E$. Assume that either $\bG$ is simple or the generic
stabilizer is semisimple.
Then $\FF(\DD(\bV(E))^{st})=\DD(\bV^{\vee}(E))^{st}_{\eps}$
(equivalently, $\FF(\SS(\bV(E))^{as})=\SS(\bV^{\vee}(E))^{as}_{\eps}$).
\end{thm}

Another result concerns the
Fourier transform of stable functions (equivalently,
antistable distributions).

\begin{thm}\label{antithm}
Let $\rho:\bG\ra\GL(\bV)$ be a nice representation
of a simply connected semisimple group $\bG$ over a $p$-adic field
$E$, such that the generic stabilizer is semisimple.
Then $\FF(\DD(\bV(E))^{as})=\DD(\bV^{\vee}(E))^{as}_{\eps}$
(equivalently, $\FF(\SS(\bV(E))^{st})=\SS(\bV^{\vee}(E))^{st}_{\eps}$).
\end{thm}

\begin{rem}
In the case of the adjoint representation the theorem \ref{mainthm}
is due to Waldspurger (for arbitrary connected reductive group $\bG$), see
\cite{Wa1}.
We conjecture that theorems \ref{mainthm} and \ref{antithm}  should hold
under a more general assumption that $E$ is a local field, and 
$\rho$ is a nice representation of a simply connected semisimple group. 
\end{rem}

Theorem \ref{mainthm} will be deduced from the more general
theorem \ref{secondmainthm}, which will be proven along
with theorem \ref{antithm} in section \ref{globalpfsec}.

The main local ingredient of these results is the analysis of
the Fourier transform of some explicit stable functions with the help of
the stationary phase principle.

\subsection{Stationary phase}
\label{statsec}

In this subsection we assume that $E$ is $p$-adic.
Let $\psi:E\ra\C^*$ be a non-trivial additive character. 
We will use the following easy version of the stationary phase principle over
$E$.

\begin{lem}\label{stationary} 
Let $\bX$ and $\bS$ be smooth varieties over $E$,
$f:\bX\times\bS\ra\A^1$ be a morphism, such that
for every $s\in\bS(\ov{E})$ the function $f_s=f|_{\bX\times\{s\}}$
on $\bX$ has finitely many non-degenerate critical points.
Let $U\subset \bX(E)$ and $P\subset \bS(E)$ be compact open subsets, 
and $\om$ a non-vanishing top-degree form on $\bX$.
For every $x_0\in\Cr(f_p)$ (where $p\in P$) we denote by 
$q_{x_0}(x-x_0)$ the quadratic form on the tangent space $T_{x_0}\bX$ 
approximating $f_p(x)-f_p(x_0)$ near $x_0$ (the Hessian of $f_p$ at $x_0$). 
Then there exists a positive constant $C$ such that 
for all $t\in (E^*)^2$ with $|t|>C$ and all $p\in P$, we have
$$
\int_U\psi(tf_p(x))|\om|=
\sum_{x_0\in\Cr(f_p)}\psi(tf_p(x_0))|t|^{-n/2}\cdot c(q_{x_0},\om_{x_0},\psi)
$$
where $n=\dim \bX$,
$c(q,\nu,\psi)=\gamma(q,\psi)\cdot |\det(q)/\nu^2|^{-1/2},$
and $\gamma(q,\psi)$ is the Weil constant associated with $q$ and $\psi$.
\end{lem}

\Pf . Let us assume first that $\bS$ is a point, so that we have
a function $f$ on $\bX$ with a finite number of
non-degenerate critical points. If $f$ has no critical points on $U$,
then we can find a finite covering $(U_i)$ of $U$ by compact open
subsets, such that on each $U_i$ there exists an analytic system of
coordinates $(x_1,\ldots,x_n)$ with $x_1=f$ and 
$\om=\la\cdot dx_1\we\ldots\we dx_n$, where $\la\in E^*$.
Furthermore, we can assume that the subsets $U_i$ are disjoint.
Therefore, the statement reduces in this case
to the vanishing of the integral
$$\int_{V}\psi(tx_1)dx_1\ldots dx_n$$
for a compact open subset $V\in E^n$ and for sufficiently large $t$,
which is clear.
Now let $c_1,\ldots,c_k$ be critical points of $f$ contained in $U$.
By Morse lemma, for each point $c_i$
there exists a neighborhood $V_i$ of $c_i$
and a system of coordinates $(x_1,\ldots,x_n)$ on $V_i$,
such that $f-f(c_i)=q(x_1,\ldots,x_n)$ for some non-degenerate
quadratic form $q$ and $\om=\la\cdot dx_1\we\ldots\we dx_n$ on $V_i$.
Let $B_i\subset V_i$ be a small ball around $c_i$ in this coordinate
system. As we have shown above,
$$\int_{U\setminus \cup_i B_i}\psi(tf)|\om|=0$$
for sufficiently large $t$. Therefore, the statement is reduced to the case
when $U$ is an open compact subgroup in a $E$-vector space $V$, $\om$ is
translation-invariant, and $f=q$ is a non-degenerate quadratic form.
Then for any $a\in E^*$ we have
$$\FF(\de_{aU})=|a|^n\cdot\vol(U)\cdot\de_{a^{-1}U^{\perp}},$$
where $U^{\perp}\subset V^{\vee}$ is the orthogonal complement to
$U$ (with respect to $\psi$).
Combining this with the equation (\ref{Weileq}) we obtain
$$\int_U\psi(a^2q)|\om|=|a|^{-n}\cdot\int_{aU}\psi(q)|\om|=
\vol(U)\cdot c(q,\om,\psi)\cdot
\int_{a^{-1}U^{\perp}}\psi(-q^{\vee})|\om|^{\vee}$$
where $|\om|^{\vee}$ is the dual measure on $V^{\vee}$.
For sufficiently large $a$ we have
$$\int_{a^{-1}U^{\perp}}\psi(-q^{\vee})|\om|^{\vee}=|a|^{-n}\cdot
\vol(U^{\perp}).$$
By involutivity of the Fourier transform, we have $\vol(U)\vol(U^{\perp})=1$,
hence we get
$$\int_U\psi(a^2q)|\om|=c(q,\om,\psi)\cdot |a|^{-n}.$$

One can deal with the case of a general family of functions
parametrized by a compact set $P\subset\bS(E)$ as follows.
The subvariety $\Cr(f)\subset\bX\times\bS$ of critical points
of $f$ in $\bX$-direction is \'etale over $\bS$. Thus,
for every point $p\in P$ the above argument works uniformly
for all $p'$ in sufficiently small neighborhood of $p$.
Now our assertion follows from compactness of $P$.
\ed

To analyze the result of the stationary phase approximation, it
is convenient to use the following lemma.

\begin{lem}\label{linearlem}
Let $f_1,\ldots,f_n$ be analytic functions on
a ball $B$ in $E^N$ centered at $0$. Assume that the
differentials at zero $d_0f_1,\ldots,d_0f_n$ are linearly independent.
Then there exists a constant $a\in E^*$, 
such that for all $t\in E^*$ with $|t|$ sufficiently large, the functions
$\psi(t f_1),\ldots,\psi(t f_n)$ on $t^{-1}a B$
are linearly independent.
\end{lem}

\Pf . Without loss of generality we can assume that $f_i(0)=0$.
Let $C\subset E^N$ be a sufficiently large open compact containing
$0$, so that the restrictions of the
functions $\psi(d_0f_1),\ldots,\psi(d_0f_n)$ to $C$ are
linearly independent. Now for sufficiently large $t$ we have
$t^{-1}C\subset B$ and $\psi(tf_i(t^{-1}x))=\psi(d_0f_i(x))$
for all $x\in C$. It remains to take $a$ such that $C\subset aB$.
\ed

\subsection{Local computation}\label{localconstrsec}

The field $E$ is still assumed to be $p$-adic.
We fix one of stabilizer subgroups $\bH\subset\bG$
and set $\bW=\bV^{\bH}$, $\bW_0=\bW\cap\bV_0$,
$\bW^{\vee}=(\bV^{\vee})^{\bH}$, $\bW^{\vee}_0=\bW^{\vee}\cap\bV^{\vee}_0$.
Let us consider the subset of $\bH$-fixed points in the variety
$\bI$:
$$\bI^{\bH}=\bI\cap(\bW_0\times\bW^{\vee}_0)=
\bI\cap(\bV_0\times\bW^{\vee}_0)$$
(the last equality follows from lemma \ref{critlem}).
Let $\bU\subset\bV_0/\bG\times\bV_0^{\vee}$ be
a non-empty $\bG$-invariant open subset 
such that the morphism
$\bff:\bI\ra\bV_0/\bG\times\bV_0^{\vee}$ is \'etale over $\bU$
(see proposition \ref{critprop}).
We have $\bI^{\bH}=\bff^{-1}(\bV_0/\bG\times\bW_0^{\vee})$. Since
$\bU$ is $\bG$-invariant, it has a non-empty intersection
with $\bV_0/\bG\times\bW_0^{\vee}$. Therefore, the open
subset $\bff^{-1}(\bU)\cap\bI^{\bH}\subset\bI^{\bH}$ is non-empty.
Note that since
$\bI_{\bH}$ is an open subset in a vector bundle over $\bW_0^{\vee}$, the
set of $E$-rational points in $\bI_{\bH}$ is dense in Zariski topology.
Therefore, there exists a point
$(x_0,x^{\vee}_0)\in \bff^{-1}(\bU)\cap\bI(E)$.
Set $\OO=\OO_{x_0}$. Then $x_0$ is a critical point of
$x^{\vee}_0|_{\OO}$ and $(\OO,x^{\vee}_0)\in\bU$.

Recall that according to lemma \ref{etale},
there exists a linear subspace $\bL\subset\bW$
such that the morphism $p:x_0+\bL\ra\bV_0/\bG$ is \'etale near $x_0$.
Let $D\subset\bL(E)$ be a small ball centered at zero,
$U^{\vee}\subset\bV_0^{\vee}$
be a compact open neighborhood of $x^{\vee}_0$. We assume that
$D$ and $U^{\vee}$ are small enough,
so that the following two conditions are satisfied:

\noindent
(i) the restriction of $p$ to $x_0+D$ is an isomorphism onto
$p(x_0+D)\subset\bV_0/\bG(E)$;

\noindent
(ii) $p(x_0+D)\times U^{\vee}\subset\bU(E)$ and
$\bff^{-1}(p(x_0+D)\times U^{\vee})$ is analytically isomorphic
to a disjoint union
of open subsets mapping identically to $p(x_0+D)\times U^{\vee}$.

Let $\OO(E)=O_1\disj\ldots\disj O_r$ be a partition of
$\OO(E)$ into $\bG(E)$-orbits. Let
us choose non-empty open compact subsets $K_i\subset O_i$, $i=1,\ldots,r$
such that $\vol(K_i)$ does not depend on $i$ (where the volume is computed
using a $\bG$-invariant top-degree form on $\OO$).
In addition we assume that for every $i=1,\ldots,r$
there exists an analytic isomorphism of $K_i$ with a ball in $E^d$,
where $d=\dim\OO$.
We set $K=K_1\disj\ldots\disj K_r\subset\OO(E)$.
Using the identification $\bG/\bH\wt{\ra}\OO:g\bH\mapsto gx_0$ we
can consider $K$ as a subset of $\bG/\bH(E)$.
By our choice of $D$,
the restriction of the map $\bG/\bH(E)\times \bL(E)\ra\bV_0(E)$
to $K\times (x_0+D)$ is an isomorphism onto 
an open compact subset $K(x_0+D)\subset\bV_0(E)$.
Now for every $t\in E^*$ with $|t|>1$ we consider a
compact open subset
$$U(t)=K(tx_0+D)=tK(x_0+t^{-1}D)\subset tK(x_0+D).$$
Note that by lemma \ref{Ginner}, stabilizers of all points in $U(t)$ are
$\bG$-inner forms of $\bH$. 
For every open compact set $C\subset\bV(E)$ we denote by $\de_{C}$
the characteristic function of $C$.
It is clear that for every $t\in E^*$ with $|t|>1$
the function $\de_{U(t)}\in\SS(\bV_0(E))$ is stable.
More generally, for every character $\kappa:A(\bH/E)\ra\C^*$
we define a function $\de^{\kappa}_{U(t)}$ supported on $U(t)$ by
\begin{equation}
\de^{\kappa}_{U(t)}(k(x_0+d))=\kappa(\inv(x_0,k))
\end{equation}
where $k\in K$, $d\in D$.

\begin{lem}\label{compsupp}
There exists a compact set $C\subset\bV^{\vee}(E)$ such
that for all $t\in E^*$ with $|t|>1$ and all $\kappa\in A(\bH/E)^D$
the support of $\FF(\de^{\kappa}_{U(t)})$
is contained in $C$.
\end{lem}

\Pf . Set $U_i(t)=K_i(tx_0+D)$, $i=1,\ldots,r$. Since
$\de^{\kappa}_{U(t)}$ is a linear combination of
$\de_{U_i(t)}$, it
suffices to prove that there exists an open compact neighborhood
of zero $B\subset\bV(E)$, such that $U_i(t)+B=U_i(t)$
for all $t\in E^*$ with $|t|>1$ and all $i$.
Let us consider the norm on $\bV(E)$ for which the unit ball
is an integer lattice in $\bV(E)$. Similarly, using an isomorphism of $K_i$
with a ball in $E^d$ and the norm on $\bL(E)$ for which $D$
is the unit ball, we get an (ultra)metric on $K_i\times (x_0+D)$ such that
$d((k,x),(k',x'))=\max(||k-k'||,||x-x'||)$.
Since the isomorphism $a:K_i\times (x_0+D)\wt{\ra} U_i=K_i(x_0+D)$ is analytic
we have $d(a^{-1}(y),a^{-1}(y'))\le c\cdot ||y-y'||$ for
some constant $c>0$, where $y,y'\in U_i$. Now for $t\in E^*$ with
$|t|>1$ the subset $t^{-1}U_i(t)=K_i(x_0+t^{-1}D)\subset U_i$ consists
of points $y\in U_i$ such that $d(a^{-1}(y),K_i\times x_0)\le |t|^{-1}$
(since $D$ is the unit ball for the norm on $\bL(E)$).
Now let $B\subset\bV(E)$ be a ball of radius $<c^{-1}$
centered at zero such that $U_i+B=U_i$. We claim that $U_i(t)+B=U_i(t)$.
Indeed, we have to check that for every $y\in t^{-1}U_i(t)$ one has  
$y+t^{-1}B\subset t^{-1}U_i(t)$. Let $y'\in y+t^{-1}B$. Then
$y'\in U_i$ since $U_i+t^{-1}B=U_i$. Also, $||y-y'||\le c^{-1}|t|^{-1}$,
hence $d(a^{-1}(y),a^{-1}(y'))\le |t|^{-1}$. By ultrametric
triangle inequality this implies that
$d(a^{-1}(y'),K_i\times x_0)\le |t|^{-1}$, i.e., $y'\in t^{-1}U_i(t)$.
\ed

\begin{prop}\label{mainloc}
(i) There exists a constant $A>0$ such that
for all $t\in (E^*)^2$ with $|t|>A$, all $x^{\vee}\in U^{\vee}$ and
$y^{\vee}\in\OO_{x^{\vee}}(E)$ one has
$$I(\FF(\de_{U(t)}))(y^{\vee})=\eps(x^{\vee},y^{\vee})\cdot
I(\FF(\de_{U(t)}))(x^{\vee}).$$

\noindent
(ii) There exists a ball $D^{\vee}\subset\bW^{\vee}(E)$
centered at zero, such that for every sufficiently small ball $D$
as above, the restriction of the function
$y\mapsto \de^{\kappa^{-1}\sign}_{\OO_y}(\FF(\de^{\kappa}_{U(t)}))$ to
$x^{\vee}_0+t^{-1}D^{\vee}$ is not identically zero provided that
$t\in (E^*)^2$ is large enough.
\end{prop}

\Pf . (i) For $x^{\vee}\in\bV_0^{\vee}(E)$ we have
$$\FF(\de_{U(t)})(x^{\vee})=\int_{U(t)}\psi(\lan x^{\vee},x\ran)dx=
|t|^{\dim \bV}\cdot\int_{t^{-1}U(t)}\psi(t\lan x^{\vee},x\ran)dx$$
where $dx$ is the Haar measure on $\bV(E)$ corresponding to a top-degree
form defined over $E$. By lemma \ref{etale}(b) 
we can rewrite this integral as follows:
$$\int_{t^{-1}U(t)}\psi(t\lan x^{\vee},x\ran)dx=
\int_{s\in x_0+t^{-1}D}\int_{k\in K}\psi(t\lan x^{\vee},ks\ran)|\om(k)|\cdot
|\nu(s)|$$
where $\om$ is a $\bG$-invariant top-degree form on $\bG/\bH$,
$\nu$ is a top-degree form on $\bL$.
Now the inner integral has form
$$\int_{x\in Ks}\psi(t\lan x^{\vee},x\ran)|\om|$$
so we can apply the stationary phase principle to compute it.
More precisely, by lemma \ref{compsupp} we know that
$\FF(\de_{U(t)})(x^{\vee})=0$ for $x^{\vee}\not\in C$, where
$C$ is a compact in $\bV^{\vee}(E)$.
Also we know that for $x^{\vee}\in p^{-1}(p(U^{\vee}))$ and
$s\in x_0+D$, the function $x^{\vee}|_{\OO_s(E)}$
has a finite number of non-degenerate critical points.
Finally, we observe that the subset $p^{-1}(p(U^{\vee}))\subset
\bV^{\vee}(E)$ is closed, hence, $C\cap p^{-1}(p(U^{\vee}))$
is compact. Thus, applying lemma \ref{stationary} we derive that
there exists a constant $A>0$ such that for $s\in x_0+D$,
$x^{\vee}\in C\cap p^{-1}(p(U^{\vee}))$, $t\in (E^*)^2$, $|t|>A$ one has
$$\int_{x\in Ks}\psi(t\lan x^{\vee},x\ran)|\om|=
|t|^{-\frac{\dim\OO}{2}}\cdot\sum_{x\in Ks\cap \Cr(x^{\vee}|_{\OO_s})}
\psi(t\lan x^{\vee},x\ran)c(B_{x,x^{\vee}},\om_x,\psi).$$
We claim that enlarging $C$ if necessary we can achieve that
the RHS is zero for 
$x^{\vee}\in p^{-1}(p(U^{\vee}))\setminus C$. Indeed, this follows
immediately from the fact that
$\bI(E)\cap (K(x_0+D)\times p^{-1}(p(U^{\vee})))$ is compact
(as a preimage of the compact set $K(x_0+D)\times p(U^{\vee})$
under the map $\bI\ra\bV_0\times \bV^{\vee}_0/\bG$). 
Thus, if $t\in (E^*)^2$ is large enough then
for all $x^{\vee}\in p^{-1}(p(U^{\vee}))$ we have
\begin{equation}\label{fourier1}
\FF(\de_{U(t)})(x^{\vee})=|t|^{\dim \bV-\frac{\dim\OO}{2}}\cdot
\int_{s\in x_0+t^{-1}D}
\sum_{x\in Ks\cap\Cr(x^{\vee}|_{\OO_s})}\psi(t\lan x^{\vee},x\ran)
c(B_{x,x^{\vee}},\om_x,\psi)|\nu(s)|.
\end{equation}
Now let us substitute $x^{\vee}$ by another
point $gx^{\vee}\in\OO_{x^{\vee}}(E)$ and make the change of
variables $x\mapsto gx$ (this makes sense since $x$ has the same stabilizer
as $x^{\vee}$). Then we obtain
$$\FF(\de_{U(t)})(gx^{\vee})=|t|^{\dim \bV-\frac{\dim\OO}{2}}\cdot
\int_{s\in x_0+t^{-1}D}
\sum_{gx\in Ks\cap\Cr(gx^{\vee}|_{\OO_s})}\psi(t\lan gx^{\vee},gx\ran)
c(B_{gx,gx^{\vee}},\om_{gx},\psi)|\nu(s)|.$$
Using lemma \ref{quadr} we can rewrite the inner sum as follows:
$$\sum_{x\in g^{-1}Ks\cap\Cr(x^{\vee}|_{\OO_s})}\psi(t\lan x^{\vee},x\ran)
\eps_{\bH}(e_{g}) c(B_{x,x^{\vee}},\om_x,\psi),$$
where $e_{g}\in H^1(E,\bH)$ is the cohomology class of the
cocycle $\sigma\mapsto g^{-1}\sigma(g)$.
Hence, we obtain
\begin{align*}
&\FF(\de_{U(t)})(gx^{\vee})=\\
&|t|^{\dim \bV-\frac{\dim\OO}{2}}\cdot
\eps_{\bH}(e_{g})\cdot\int_{s\in x_0+t^{-1}D}
\sum_{x\in g^{-1}Ks\cap\Cr(x^{\vee}|_{\OO_s})}\psi(t\lan x^{\vee},x\ran)
c(B_{x,x^{\vee}},\om_x,\psi)|\nu(s)|.
\end{align*}
Now to calculate $I(\FF(\de_{U(t)}))$ at $gx^{\vee}$, 
we have to replace $g$ by $g_1g$ in the above formula
where $g_1\in\bG(E)/(g\bH g^{-1})(E)$ and integrate over
$g_1$. We get
\begin{align*}
&I(\FF(\de_{U(t)}))(gx^{\vee})=
|t|^{\dim \bV-\frac{\dim\OO}{2}}\cdot
\eps_{\bH}(e_{g})\times\\
&\int_{s\in x_0+t^{-1}D}
\int_{g_1\in\bG(E)/(g\bH g^{-1})(E)}|\om(g_1)|\cdot
\sum_{x\in g^{-1}g_1^{-1}Ks\cap\Cr(x^{\vee}|_{\OO_s})}
\psi(t\lan x^{\vee},x\ran)
c(B_{x,x^{\vee}},\om_x,\psi)|\nu(s)|=\\
&|t|^{\dim \bV-\frac{\dim\OO}{2}}\cdot
\eps_{\bH}(e_{g})\times\\
&\int_{s\in x_0+t^{-1}D}
\sum_{x\in\Cr(x^{\vee}|_{\OO_s})}
\vol(\bG(E)gx\cap Ks)\psi(t\lan x^{\vee},x\ran)
c(B_{x,x^{\vee}},\om_x,\psi)|\nu(s)|
\end{align*}
where volumes are computed using $\om$.
By our choice of $K$ we have
$\vol(\bG(E)gx\cap Ks)=\vol(\bG(E)x_0\cap K)$.
Hence, we conclude that
\begin{equation}\label{fouriercomp}
\begin{array}{l}
I(\FF(\de_{U(t)}))(gx^{\vee})=
|t|^{\dim \bV-\frac{\dim\OO}{2}}\cdot\vol(\bG(E)x_0\cap K)\cdot
\eps_{\bH}(e_{g})\times\\
\int_{s\in x_0+t^{-1}D}\sum_{x\in\Cr(x^{\vee}|_{\OO_s})}
\psi(t\lan x^{\vee},x\ran)c(B_{x,x^{\vee}},\om_x,\psi)|\nu(s)|.
\end{array}
\end{equation}
Since $\eps_{\bH}(e_{g})=\eps(x^{\vee},gx^{\vee})$
this finishes the proof of the part (i) of the proposition.

\noindent
(ii) Arguing as above, we obtain that
\begin{equation}\label{fouriercomp2}
\begin{array}{l}
I(\FF(\de^{\kappa}_{U(t)}))(gx^{\vee})=
|t|^{\dim \bV-\frac{\dim\OO}{2}}\cdot\vol(\bG(E)x_0\cap K)\cdot
\eps_{\bH}(e_{g})\times\\
\int_{s\in x_0+t^{-1}D}\sum_{x\in\Cr(x^{\vee}|_{\OO_s})}
\kappa(\inv(s,gx))\cdot
\psi(t\lan x^{\vee},x\ran)c(B_{x,x^{\vee}},\om_x,\psi)|\nu(s)|.
\end{array}
\end{equation}
for sufficiently large $t\in (E^*)^2$, where $x^{\vee},gx^{\vee}\in
p^{-1}(p(U^{\vee}))$.

By our choice of $D$ and $U^{\vee}$
there exists a collection of analytic maps
$$x_i:(x_0+D)\times U^{\vee}\ra\bV(E), i=1,\ldots,n,$$
such that for every $s\in x_0+D$
the points $x_1(s,x^{\vee}),\ldots,x_n(s,x^{\vee})$ are
disjoint and constitute the set of critical points
of $x^{\vee}|_{\OO_s}$. Let us set $x_i=x_i(x_0,x_0^{\vee})$.
Renumbering these maps if necessary we can assume that $x_1=x_0$.
We claim that the differentials at $x^{\vee}_0$ of
the functions $x^{\vee}\mapsto\lan x^{\vee},x_i(x_0,x^{\vee})\ran$,
where $x^{\vee}\in U^{\vee}\cap\bW^{\vee}(E)$, $i=1,\ldots,n$, are linearly
independent. Note that by lemma \ref{critlem}(a)
for $x^{\vee}\in U^{\vee}\cap\bW^{\vee}(E)$
we have $x_i(x_0,x^{\vee})\in\bW(E)$. Let $L_i\in\Hom(\bW^{\vee},\bW)$
be the differential of $x_i(x_0,?)|_{U^{\vee}\cap\bW^{\vee}(E)}$ at
$x_0^{\vee}$. Differentiating the condition
$x_i(x_0,x^{\vee})\in\OO_{x_0}$,
we get that $L_i(\bW^{\vee})\subset\gg x_0$.
Now the differential at $x^{\vee}_0$
of the function $\lan x^{\vee},x_i(x_0,x^{\vee})\ran$
on $U^{\vee}\cap\bW^{\vee}$ is the functional
$$w^{\vee}\mapsto \lan w^{\vee},x_i(x_0,x^{\vee}_0)\ran+
\lan x^{\vee}_0, L_i(w^{\vee})\ran=\lan w^{\vee},x_i(x_0,x^{\vee}_0)\ran$$
(here we used the equality $\lan x^{\vee}_0,\gg x_0\ran=0$).
Since all the points $x_i(x_0,x^{\vee}_0)\in\bW(E)$ are distinct, 
this proves our claim.
By lemma \ref{linearlem} this implies that there exists a ball
$D^{\vee}\subset\bW^{\vee}(E)$ centered at zero, such that
for sufficiently large $t$ the functions
\begin{equation}\label{functions}
x^{\vee}\mapsto \psi(t\lan x^{\vee},x_i(x_0,x^{\vee})\ran),\ i=1,\ldots n,
\end{equation}
are linearly independent on $x^{\vee}_0+t^{-1}D^{\vee}$.

Since the functions $x_i$ are analytic,
we can choose $D$ sufficiently small so that
$$\psi(t\lan x^{\vee},x_i(s,x^{\vee})\ran)=
\psi(t\lan x^{\vee},x_i(x_0,x^{\vee})\ran),$$
for $s\in x_0+t^{-1}D$, $x^{\vee}\in U^{\vee}$, $i=1,\ldots,n$.
On the other hand, if $t$ is large enough then
setting 
$c_i(s,x^{\vee})= 
c(B_{x_i(s,x^{\vee}),x^{\vee}},\om_{x_i(s,x^{\vee})},\psi)$ for
$i=1,\ldots,n$ , we get
$$c_i(s,x^{\vee})=c_i(x_0,x_0^{\vee}),$$
$$\kappa(\inv(s,gx_i(s,x^{\vee})))=\kappa(\inv(x_0,gx_i))$$
for $s\in x_0+t^{-1}D$, $x^{\vee}\in x^{\vee}_0+t^{-1}D^{\vee}$,
$i=1,\ldots,n$.
Finally, we have
$$\kappa(\inv(x_0,gx_i))=\kappa(\inv(x_0,x_i))\cdot\kappa(\tau_i(\inv(e_g)))$$
where $\tau_i:A(\bH/E)\ra A(\bH/E)$ is the automorphism
induced by the action of an element $g_i\in N(\bH)/\bH$ such that 
$g_ix_i=x_0$. Note that since by our assumption $x_1=x_0$, we have
$\tau_1=\id$.
Thus, applying the formula (\ref{fouriercomp2})
to $x^{\vee}\in x^{\vee}_0+t^{-1}D^{\vee}$ and large enough
$t\in (E^*)^2$, and using Lemma \ref{signinv}, we get
\begin{align*}
&I(\FF(\de_{U(t)}))(gx^{\vee})=
|t|^{\dim \bV-\frac{\dim\OO}{2}}\cdot
\vol(\bG(E)x_0\cap K)\times\\
&\int_{s\in x_0+t^{-1}D}\sum_{i=1}^n
\kappa(\inv(x_0,x_i))((\kappa\circ\tau_i)\sign)(\inv(e_g))
\psi(t\lan x^{\vee},x_i(s,x^{\vee})\ran)
c_i(x_0,x_0^{\vee})|\nu(s)|=\\
&|t|^{\dim \bV-\frac{\dim\OO}{2}}\cdot
\vol(\bG(E)x_0\cap K)\vol(t^{-1}D)\times\\
&\sum_{i=1}^n
\kappa(\inv(x_0,x_i))((\kappa\circ\tau_i)\sign)(\inv(e_g))
c_i(x_0,x_0^{\vee})\psi(t\lan x^{\vee},x_i(x_0,x^{\vee})\ran).
\end{align*}
Hence, we have
\begin{align*}
&\delta^{\kappa^{-1}\sign}_{\OO_{x^{\vee}}}(\FF(\de_{U(t)}))=
c(t)\cdot \sum_{gx^{\vee}\in\bG(E)\backslash\OO_{x^{\vee}}(E)}\\
&\sum_{i=1}^n
\kappa(\inv(x_0,x_i))\frac{\kappa\circ\tau_i}{\kappa}(\inv(e_g))
c_i(x_0,x_0^{\vee})\psi(t\lan x^{\vee},x_i(x_0,x^{\vee})\ran),
\end{align*}
where $c(t)\in\C^*$ is a constant depending on $t$.
Interchaning two summations and noting that $\inv(e_g)$ runs
through the entire group $A=A(\OO_{x^{\vee}})=A(\OO_{x^{\vee}_0})$, we
obtain
\begin{equation}
\delta^{\kappa^{-1}\sign}_{\OO_{x^{\vee}}}(\FF(\de_{U(t)}))=
c(t)|A|\cdot\sum_{i:\kappa\circ\tau_i=\kappa}\kappa(\inv(x_0,x_i))
c_i(x_0,x_0^{\vee})\psi(t\lan x^{\vee},x_i(x_0,x^{\vee})\ran)
\end{equation}
for $x^{\vee}\in x^{\vee}_0+t^{-1}D^{\vee}$ and large enough
$t\in (E^*)^2$.
Now our assertion follows from
linear independence of the functions (\ref{functions}) on
$x^{\vee}_0+t^{-1}D^{\vee}$ and the fact that $\tau_1=\id$.
\ed

\subsection{Inner forms}\label{innersec}

Let $\bZ\subset \bG$ be a central subgroup defined over $E$,
$\a\in H^1(E,\bG/\bZ)$
be a cohomology class. Then $\a$ defines an inner form $\bG'$ of $\bG$.
Note that since $\bG$ is simply connected, the homomorphism
$$d:H^1(E,\bG/\bZ)\ra H^2(E,\bZ)$$
is an isomorphism when $E$ is $p$-adic,
so in this case $\a$ is uniquely determined by $d(\a)\in H^2(E,\bZ)$.
Assume that $\rho(\bZ)=1$. Then we can twist $\rho$ by $\a$ to get
a nice representation $\rho':\bG'\ra\Aut(\bV')$.
Let $\bH\subset\bG$ (resp. $\bH'\subset\bG'$) be a generic stabilizer
subgroup, $\bV_0\subset \bV$ (resp. $\bV'_0$)
be the open subset consisting of points with
the stabilizer conjugated to $\bH$ (resp. $\bH'$) over $\ov{E}$.
The isomorphism $i:\bV\ra\bV'$ defined over $\ov{E}$ induces a bijection
between $\bG$-orbits on $\bV$ defined over $E$ and $\bG'$-orbits on $\bV'$
defined over $E$. Let $\OO\subset\bV_0$ be a $\bG$-orbit,
$\OO'=i(\OO)\subset\bV'_0$ be the corresponding $\bG'$-orbit.
If $\om$ is a $\bG$-invariant non-zero top degree form on $\OO$ defined
over $E$ then $i_*\om$ is a $\bG'$-invariant top degree form on $\OO'$,
also defined over $E$.

\noindent {\bf Definition}.
For a pair of functions $\phi\in\SS(\bV(E))$, $\phi'\in\SS(\bV'(E))$ we say
that $\phi\sim \phi'$ (resp. $\phi\sim_{\eps} \phi'$)
if for every $\bG$-orbit $\OO\subset\bV_0$ defined
over $E$ one has $\delta_{\OO,\om}(\phi)=\delta_{\OO',i_*\om}(\phi')$
(resp. $\delta^{\eps}_{\OO,\om}(\phi)=\delta^{\eps}_{\OO',i_*\om}(\phi')$).
For a pair of distributions $\delta\in\DD(\bV(E))$, $\delta'\in\SS(\bV'(E))$
we say that $\delta\sim\delta'$ (resp. $\delta\sim_{\eps}\delta'$)
if for every pair of functions $\phi\in\SS(\bV(E))$, $\phi'\in\SS(\bV'(E))$
such that $\phi\sim \phi'$ (resp. $\phi\sim_{\eps} \phi'$) one has
$\delta(\phi)=\delta(\phi')$.

\vspace{2mm}

Note that we have $\phi\sim 0$ (resp. $\phi\sim_{\eps} 0$) if and only if
$\phi$ is antistable (resp. $\eps$-antistable). Therefore, if
$\delta\sim\delta'$ (resp. $\delta\sim_{\eps}\delta'$)
then distributions $\delta$ and $\delta'$ are necessarily
stable (resp. $\eps$-stable).
Also, by the definition we have $\delta_{\OO,\om}\sim\delta_{\OO',i_*\om}$
(resp. $\delta^{\eps}_{\OO,\om}\sim_{\eps}\delta^{\eps}_{\OO',i_*\om}$).

To formulate our result on inner forms we have to introduce
certain sign associated with $\bG$ and $\bG'$.
Let $Q$ (resp. $Q'$) be the Killing form on $\gg$ (resp. $\gg'$). We set
$$\kappa(\bG,\bG')=\eps(Q,Q').$$
Note that $Q$ and $Q'$ have the same determinant modulo squares.
Thus, in the case of $p$-adic $E$
the difference between equivalence classes of these
quadratic forms is measured by the sign $\kappa(\bG,\bG')$.

\begin{thm}\label{secondmainthm}
Let $\rho:\bG\ra\GL(\bV)$ be a nice representation
of a simply connected semisimple group $\bG$ over a $p$-adic field
$E$, $\bZ\subset\bG$ be a central subgroup acting trivially on $\bV$.
Let $(\bG',\bV')$ be a twist
of $(\bG,\bV)$ by a class $\a\in H^1(E,\bG/\bZ)$.
Assume that either $\bG$ is simple or the generic stabilizer 
is semisimple.
Then for every pair of functions $\phi\in\SS(\bV(E))$, $\phi'\in\SS(\bV(E))$
one has $\phi\sim \phi'$ if and only if
$\FF(\phi)\sim_{\eps}\kappa(\bG,\bG')\FF(\phi')$.
The same result holds for distributions.
\end{thm}

The proof will be given in \ref{globalpfsec}.
The remainder of this section consists of various local ingredients
of the proof.

Theorem \ref{mainthm} is an immediate consequence
of theorem \ref{secondmainthm}. Indeed, we can take $(\bG',\bV')=(\bG,\bV)$
and $\phi'=0$. Then the condition $\phi\sim 0$ means that $\phi$ is antistable,
while the condition $\FF(\phi)\sim_{\eps} 0$ means that $\FF(\phi)$ is
$\eps$-antistable.

\begin{lem}\label{nonempty}
(a) Let $\OO\subset\bV_0$ be an orbit, $x\in\OO(E)$
be an $E$-point on $\OO$.
Then $\OO'(E)$ is non-empty if and only if the class $\a$ belongs
to the image of the map $H^1(E,\bH_x/\bZ)\ra H^1(E,\bG/\bZ)$.

\noindent (b) For every $\phi\in\SS(\bV'(E))$ the functions
$x\mapsto \delta_{\OO'_x,i_*\om}(\phi)$ and 
$x\mapsto\delta^{\eps}_{\OO'_x,i_*\om}(\phi)$
on $\bV_0(E)$ are locally constant.
\end{lem}

\Pf . The proof of (a) is straightforward. To prove (b) let us fix
a point $x\in\bV_0(E)$. Set $\bH=\bH_x$, $\bW=\bV^{\bH}$,
$\bW_0=\bW\cap\bV_0$. The morphism $a:\bG/\bH\times\bW_0\ra\bV_0$ is
smooth, hence $U=a(\bG/\bH(E)\times\bW_0(E))$ is an open
subset in $\bV_0(E)$.
Assume first that $\a$ does not belong to the
image of the map $H^1(E,\bH/\bZ)\ra H^1(E,\bG/\bZ)$. Then
for every point $y\in U$ we have $\OO'_y(E)=\emptyset$. Thus,
we can assume that $\a$ belongs to the image of this map.
Then there exists an isomorphism $i:(\bG,\bV)\ra(\bG',\bV')$
over $\ov{E}$, such that the subgroup $i({\bH})\subset\bG'$ and
the morphism $i|_{\bW}$ are defined over $E$.
Therefore, the functions $y\mapsto \delta_{\OO'_y,i_*\om}(\phi)$ 
and $y\mapsto \delta^{\eps}_{\OO'_y,i_*\om}(\phi)$ on $\bW_0(E)$
is locally constant. Let $U_0\subset \bW_0(E)$ be a neighborhood of
$x$ on which it is constant. Then $a(\bG/\bH(E)\times U_0)$ is an open
neighborhood of $x$ in $\bV_0(E)$ on which a similar function is constant.
\ed

\begin{lem}\label{anisotrlem}
Assume that $E$ is $p$-adic, $\bG$ is simple,
the generic stabilizer for $\rho$ is commutative and
all irreducible components
of $\rho_{\ov{E}}$ are defined over $E$. Then there exists a point
$x\in\bV_0(E)$ such that $\bH_x$ is an anisotropic.
\end{lem}

\Pf .
In the case when $\rho$ is the adjoint representation,
we can use the existence of an anisotropic maximal torus $\bH_x$ in $\bG$
defined over $E$ (see e.g., Theorem VI.21 of \cite{PR}). 
It turns out that in all other cases the action of a sufficiently
big subgroup of $\bG$ on $\bV$ essentially reduces to the adjoint
representation. Here is a more precise statement.

\noindent
{\it Claim}. For every $x\in\bV_0(E)$
there exists a semisimple subgroup $\bG'\subset\bG$ defined over $E$
and a $\bG'$-invariant decomposition $\bV=\bV'\oplus\bV''$, such that
the following two conditions hold:

\noindent
(i) the representation of $\bG'$ on $\bV''$ is equivalent to the 
adjoint representation of $\bG'$;

\noindent
(ii) let $x=x'+x''$ where $x'\in\bV'$, $x''\in\bV''$, then $\bG' x'=x'$
and $\bH_x$ coincides with the stabilizer of $x''$ in $\bG'$
(so by (i), $\bH_x$ is a maximal torus in $\bG'$).

\vspace{1mm}

Our statement can be deduced from this claim 
as follows. The subset $\bV_0\cap (x'+\bV'')\subset x'+\bV''$ is non-empty 
and Zariski open.
Therefore, we can choose $\wt{x}''\in\bV''(E)$ 
such that $\wt{x}=x'+\wt{x}''\in\bV_0(E)$
and the stabilizer of $\wt{x}''$ in $\bG'$ is anisotropic. Then
$\bH_{\wt{x}}$ contains an anisotropic maximal torus in
$\bG'$. Since $\wt{x}$ is contained in $\bV_0$, 
this inclusion is in fact an equality.  
The proof of the Claim follows from Elashvili's classification
of representations of simple groups with
generic stabilizers of positive dimension (see \cite{El1}).
Here are the cases relevant for our situation:

\noindent
(i) $\bG_{\ov{E}}=\SL(W)$, $\bV_{\ov{E}}=\bV'\oplus\bV''$, where
$\bV'$ and $\bV''$ are irreducible representations of $\SL(W)$
in $S^2W$ and $\We^2 W$ or in $S^2 W^{\vee}$ and $\We^2 W$ (or dual to
these).
The stabiliser $\bG'$ of a generic point in $\bV'$ is the special orthogonal
group $\SO(W)$. It is well-known that the representation of $\SO(W)$
in $\We^2 W$ is isomorphic to the adjoint representation.

\noindent
(ii) $\bG_{\ov{E}}=\SL(W)$, where $\dim W=4$,
$\bV_{\ov{E}}$ is the direct sum of $4$ copies of $\We^2 W$.
The stabilizer in $\SL(W)$ of a generic point in $(\We^2 W)^2$
is conjugate to the subgroup $\SL_2\times\SL_2$ corresponding
to a decomposition $W=W_1\oplus W_2$, where $\dim W_1=\dim W_2=2$
(see Table 1 of \cite{El1}).
The generic stabilizer in $(\We^2 W)^3$ is the subgroup
$\SL_2=\{(g,g^{-1}), g\in\SL_2\}$ in $\SL_2\times\SL_2$ (which corresponds
to choosing an isomorphism $W_1\simeq W_2$).
Decomposing the $4$-th factor as a representation of $\SL(W_1)\times\SL(W_2)$:
$$\sideset{}{^2}{\We} W=1\oplus 1\oplus W_1\otimes W_2,$$
we see that the action of the above subgroup on $\We^2 W$
is equivalent to the direct sum of the adjoint representation of $\SL_2$
with $3$ trivial representations. Thus, the above claim holds if we
take $\bG'$ to be the stabilizer of the first three components of $x$,
$\bV'$ to be the sum of $(\We^2 W)^3$ and of $3$ trivial
representations.

\noindent
(iii) $\bG_{\ov{E}}=\SL(W)$, where $\dim W=6$,
$\bV_{\ov{E}}$ is the direct sum of $2$ copies of $\We^3 W$.
The stabiliser in $\SL(W)$ of a generic point in $\We^3 W$
is conjugate to the subgroup $\SL_3\times\SL_3$ corresponding
to a decomposition $W=W_1\oplus W_2$, where $\dim W_1=\dim W_2=3$
(see Table 1 of \cite{El1}).
We can decompose the second copy of $\We^3 W$ as a representation of
$\SL(W_1)\times\SL(W_2)$:
$$\sideset{}{^3}{\We} W=1\oplus 1\oplus W_1^{\vee}\otimes W_2\oplus
W_2^{\vee}\otimes W_1.$$
Note that the last two factors are non-isomorphic non-trivial irreducible
representations of $\SL(W_1)\times\SL(W_2)$, so for any $E$-form of the pair
$(\SL(W_1)\times\SL(W_2),\We^3 W)$ the similar decomposition takes place.
Now the stabilizer in $\SL(W_1)\times\SL(W_2)$ of a generic point
in one of the non-trivial
factors is isomorphic to $\SL_3$, and the action of this subgroup on the
second non-trivial factor is the sum of the adjoint representation and
of the trivial representation. Using this we can easily construct $\bV'$ and
$\bV''$ (with $\bG'$ being the form of $\SL_3$).

\noindent
(iv) $\bG_{\ov{E}}=\SL(W)$, where $\dim W=8$,
$\bV_{\ov{E}}$ is either $\We^3 W\oplus W$ or $\We^3 W\oplus W^{\vee}$.
The stabilizer in $\SL(W)$ of a generic point in $\We^3 W$ is isomorphic
to $\SL_3$. The embedding of $\SL_3$ in $\SL(W)$ corresponds to the 
identification of $W$ with the adjoint representation of $\SL_3$ (see Table 1
of \cite{El1}). Therefore, we can take $\bG'$ to be 
the stabilizer of the first component of $x$ with respect to the above
decomposition of $\bV$. 
\ed

\begin{rem} Note that if $\bG$ is split over $E$, then
all irreducible components
of $\rho_{\ov{E}}$ are defined over $E$.
\end{rem}

\begin{lem}\label{maininner}
Assume that $E$ is $p$-adic and
 that either $\bH$ is semisimple, or
$\bG$ is simple and all irreducible components
of $\rho_{\ov{E}}$ are defined over $E$.
Then there exists a point $x\in\bV_0(E)$, such that
the map $H^1(E,\bH_x/\bZ)\ra H^1(E,\bG/\bZ)$ is surjective.
\end{lem}

\Pf . Note that the natural map $H^1(E,\bG/\bZ)\ra H^2(E,\bZ)$
is an isomorphism (since $E$ is $p$-adic).
Therefore, we have to prove that the natural map
$$H^1(E,\bH_x/\bZ)\ra H^2(E,\bZ)$$
is surjective for some $x\in\bV_0(E)$.

Assume first that $\bH$ is semisimple. Then we claim that this
surjectivity holds for every point $x\in\bV_0(E)$. Indeed,
let $\wt{\bH}\ra \bH$ be the universal covering of $\bH$,
$\wt{\bZ}\subset\wt{\bH}$ be the preimage of $\bZ\subset\bH$.
Then we have an isomorphism
$$H^1(E,\bH/\bZ)=H^1(E,\wt{\bH}/\wt{\bZ})\wt{\ra} H^2(E,\wt{\bZ}).$$
Since the cohomological dimension of $E$ is equal to $2$, the map
$H^2(E,\wt{\bZ})\ra H^2(E,\bZ)$ is surjective, which finishes the proof
in this case (the same argument can be applied to any $\bH_x$).

Now let us assume that $\bH$ is commutative. Then by lemma \ref{anisotrlem}
there exists a point $x\in\bV_0(E)$ such that $\bH_x$ is an anisotropic
torus.
By Tate-Nakayama duality, for such a torus we have $H^2(E,\bH_x)=0$,
which implies the surjectivity we want.
\ed

In the remainder of this section
we will keep the assumptions of lemma \ref{maininner}.
Furthermore, we choose a point $x_0\in\bV_0(E)$ as in this lemma
and set $\bH=\bH_{x_0}$, $\bW=\bV^{\bH}$, $\bW_0=\bW\cap\bV_0$, etc.
We also fix a cohomology class $\a_0\in H^1(E,\bH/\bZ)$
mapping to $\a$. This allows us to choose
an isomorphism $i:(\bG,\bV)\ra (\bG',\bV')$
over $\ov{E}$, such that the subgroup $\bH'=i(\bH)\subset\bG'$
and the morphism $i|_{\bW}$ are defined over $E$.
We denote $\bW'=(\bV')^{\bH'}=i(\bW)$.  Also we denote by
$i^{\vee}:\bV^{\vee}\ra(\bV')^{\vee}$ the $\ov{E}$-isomoprhism induced
by $i$. Note that $i^{\vee}|_{\bW^{\vee}}$ is also defined over $E$.

The following lemma (generalizing lemma \ref{quadr})
computes the sign that in the $p$-adic case measures the difference
between the quadratic form $B_{x,x^{\vee}}$ introduced in
\ref{spinsec} and $B_{i(x),i^{\vee}(x^{\vee})}$.
Recall that for $x\in\OO$, $x'\in\OO'$ we denote
$\eps(x,x')=\eps(Q_x,Q'_{x'})$.

\begin{lem}\label{signlem}
Let $x\in\bW_0$ be a critical point of $x^{\vee}|_{\OO_x}$,
where $x^{\vee}\in\bW_0^{\vee}$.
Let $\om$ be a non-zero $\bG$-invariant top-degree form on $\OO_x$
defined over $E$. 

\noindent
(a) One has $\det(B_{i(x),i^{\vee}(x^{\vee})})/(i_*\om)_{i(x)}^2=
\det(B_{x,x^{\vee}})/\om_x^2$.

\noindent
(b) The quadratic forms $B_{i(x),i^{\vee}(x^{\vee})}$ and $B_{x,x^{\vee}}$
have the same determinant modulo squares. Their relative Hasse-Witt invariant
is given by
$$\eps(B_{i(x),i^{\vee}(x^{\vee})},B_{x,x^{\vee}})=
\kappa(\bG,\bG')\eps(x,i(x))=\kappa(\bG,\bG')\eps(x^{\vee},i(x^{\vee})).$$
\end{lem}

\Pf. The proof of (a) is straightforward. To prove (b) we
note that there is a natural homomorphism
$\iota:\bH/\bZ\ra\SO(T_x,B_{x,x^{\vee}})$ which induces a map
$$\iota_*:H^1(E,\bH/\bZ)\ra H^1(E,\SO(T_x, B_{x,x^{\vee}})).$$
It is easy to see that the quadratic form $B_{i(x),i^{\vee}(x^{\vee})}$
is equivalent to the twist of $B_{x,x^{\vee}}$ by $\iota_*(\a_0)$.
In particular, these forms have the same determinant modulo
squares. Let
$$\delta:H^1(E,\SO(T_x,B_{x,x^{\vee}}))\ra H^2(E,\{\pm 1\})\simeq\{\pm 1\}$$
be the map induced by the spin-covering of $\SO(T_x, B_{x,x^{\vee}})$.
Then the Hasse-Witt invariants of $B_{v,v^{\vee}}$ and
of its twist by $\iota_*(\a_0)$ differ by $\delta(\iota_*(\a_0))$.
It remains to prove that
\begin{equation}\label{coh1}
\delta(\iota_*(\a_0))=\kappa(\bG,\bG')\eps(x,i(x)).
\end{equation}
Recall that by proposition \ref{spinprop} the homomorphism
$\iota$ lifts to a homomorphism
$$\wt{\bH}\ra\Spin(T_x, B_{x,x^{\vee}})$$
where $\wt{\bH}\ra\bH$ is the pull-back of the spin-covering of
$\SO(\hh,Q|_{\hh})$. Let $\wt{\bZ}\subset\wt{\bH}$ be the preimage of
$\bZ\subset\bH$. Then we have an induced homomorphism
$$\chi:\wt{\bZ}\ra\{\pm 1\}.$$
On the other hand, we have a natural map of cohomologies
$$\wt{d}:H^1(E,\bH/\bZ)=H^1(E,\wt{\bH}/\wt{\bZ})\ra H^2(E,\wt{\bZ}).$$
Now it is easy to see that
$$\delta\circ\iota_*=\chi_*\circ\wt{d}$$
where $\chi_*:H^2(E,\wt{\bZ})\ra H^2(E,\{\pm 1\})$ is
the homomorphism induced by $\chi$. Therefore, we have
\begin{equation}\label{coh2}
\delta(\iota_*(\a_0))=\chi_*(\wt{d}(\a_0))
\end{equation}
On the other hand, we have the natural commutative diagram
\begin{equation}
\begin{array}{ccc}
\bG & \lrar{} & \Spin(\gg,Q)\\
\ldar{} & & \ldar{}\\
\bG/\bZ & \lrar{} & \SO(\gg,Q)
\end{array}
\end{equation}
so we get a homomorphism $\chi'_{\bG}:\bZ\ra\{\pm 1\}$ such that
\begin{equation}\label{coh2h}
\kappa(\bG,\bG')=\chi'_{\bG,*}(d(\alpha)).
\end{equation}
Composing it
with the natural projection $\wt{\bZ}\ra\bZ$ we get a homomorphism
$\chi_{\bG}:\wt{\bZ}\ra\{\pm 1\}$ such that
\begin{equation}\label{coh3}
\kappa(\bG,\bG')=\chi_{\bG,*}(\wt{d}(\alpha_0)).
\end{equation}
Also, by the definition of $\wt{\bH}$ we have the following commutative
diagram
\begin{equation}
\begin{array}{ccc}
\wt{\bH} & \lrar{} & \Spin(\hh,Q|_{\hh})\\
\ldar{} & & \ldar{}\\
\wt{\bH}/\wt{\bZ} & \lrar{} & \SO(\hh,Q|_{\hh})
\end{array}
\end{equation}
which gives the homomorphism $\chi_{\bH}:\wt{\bZ}\ra\{\pm 1\}$ such that
\begin{equation}\label{coh4}
\eps(x,i(x))=\chi_{\bH,*}(\wt{d}(\alpha_0)).
\end{equation}
It remains to prove the equality
$$\chi\cdot\chi_{\bG}\cdot\chi_{\bH}=1$$
of homomorphisms from $\wt{\bZ}$ to $\{\pm 1\}$.
Indeed, then the equation (\ref{coh1}) would follow from
(\ref{coh2}), (\ref{coh3}) and (\ref{coh4}).
To prove the equality of algebraic homomorphisms we can pass
to the algebraic closure of $E$. Hence, it suffices to do this
over $\C$. Then we can proceed similarly to the proof of proposition
\ref{spinprop}. Namely, as in that proof we deduce that
$\chi$ coincides with the homomorphism induced by
the homomorphism $\wt{\bH}\ra\Spin(\hh^{\perp},Q|_{\hh^{\perp}})$ (the
existence of such a homomorphism also follows from the proof of proposition
\ref{spinprop}).
Now our statement follows from the commutativity
of the diagram
\begin{equation}
\begin{array}{ccc}
\wt{\bH} & \lrar{} & \Spin(\hh,Q|_{\hh})\times 
\Spin(\hh^{\perp},Q|_{\hh^{\perp}})\\
\ldar{} & &\ldar{}\\
\bG & \lrar{} & \Spin(\gg,Q)
\end{array}
\end{equation}
where the left vertical arrow factors through $\bH$.
\ed

\begin{rem} It is well-known that the spinor representation of
$\bG$ is a multiple of the irreducible representation $V_{\rho}$
corresponding to the half-sum of positive roots $\rho$ (see \cite{Kos}).
It follows that the character $\chi'_{\bG}$ concides with the restriction
to $\bZ$ of $\rho$ considered as a character of the maximal
torus of $\bG$. Together with the equation (\ref{coh2h}) 
this implies that the sign $\kappa(\bG,\bG')$ coincides
with the sign introduced by Kottwitz in \cite{K0}.
\end{rem}

The following proposition shows how to construct pairs
of functions $\phi\in\SS(\bV_0(E))$, $\phi'\in\SS(\bV'_0(E))$
with $\phi\sim \phi'$ and with
some control over the Fourier transforms of $\phi$ and $\phi'$.

\begin{prop}\label{innerprop}
Keep the assumptions of lemma \ref{maininner}. Then there exist
$\phi\in\SS(\bV_0(E))$ and $\phi'\in\SS(\bV'_0(E))$
such that $\phi\sim \phi'$, while
\begin{equation}\label{innerstat}
\delta^{\eps}_{\OO_{x^{\vee}},\om}(\FF(\phi))=
\kappa(\bG,\bG')\delta^{\eps}_{\OO'_{x^{\vee}},i_*\om}(\FF(\phi'))\neq 0
\end{equation}
for some $x^{\vee}\in\bV^{\vee}_0(E)$.
\end{prop}

\Pf . We use the notation of section \ref{localconstrsec}.
Let us set
$$\phi=\frac{1}{\vol_{\om}(\bG(E)x_0\cap K)}\cdot\delta_{U(t)}$$
for $t\in (E^*)^2$ sufficiently large.
Recall that $U(t)=K(tx_0+D)\subset\bV_0(E)$. Let us choose
some compact open subset
$K'\subset\OO_{i(x_0)}(E)$ intersecting $\bG'(E)$-orbits at subsets of
equal volumes, and let $U'(t)=K'(ti(x_0)+i(D))$ be the corresponding
compact open subset of $\bV'(E)$. Now let us set
$$\phi'=\frac{1}{\vol_{i_*\om}(\bG'(E)i(x_0)\cap K')}\cdot\delta_{U'(t)}.$$
Then clearly we have $\phi\sim \phi'$. 
On the other hand, using the formula
(\ref{fouriercomp}) we obtain that for $x^{\vee}\in U^{\vee}\cap
\bW^{\vee}_0(E)$
one has
\begin{align*}
&\delta^{\eps}_{\OO_{x^{\vee}}}(\FF(\phi))=
|t|^{\dim \bV-\frac{\dim\OO}{2}}\cdot
|H^1(E,\bH)|\cdot\eps(x^{\vee})\times\\
&\int_{s\in x_0+t^{-1}D}\sum_{x\in\Cr(x^{\vee}|_{\OO_s})}
\psi(t\lan x^{\vee},x\ran)c(B_{x,x^{\vee}},\om_x,\psi)|\nu(s)|,
\end{align*}
\begin{align*}
&\delta^{\eps}_{\OO_{i(x^{\vee})}}(\FF(\phi'))=
|t|^{\dim \bV-\frac{\dim\OO}{2}}\cdot
|H^1(E,\bH')|\cdot\eps(i(x^{\vee}))\times\\
&\int_{s\in x_0+t^{-1}D}\sum_{x'\in\Cr(i(x^{\vee})|_{\OO_{i(s)}})}
\psi(t\lan i(x^{\vee}),x'\ran)c(B_{x',i(x^{\vee})},\om_{x'},\psi)|\nu(s)|.
\end{align*}
where $\bH'=i(\bH)$. Since $\bH'$ is an inner form of $\bH$, 
we have $|H^1(E,\bH)|=|H^1(E,\bH')|$. On the other hand,
the isomorphism $i|_{\bW}$ sends $\Cr(x^{\vee}|_{\OO_s})$ to
$\Cr(i(x^{\vee})|_{\OO_{i(s)}})$. Hence, applying
lemma \ref{signlem} we get the equality (\ref{innerstat}).
It remains to note that 
according to proposition \ref{mainloc}(ii), 
we will also have $\delta^{\eps}_{\OO_{x^{\vee}},\om}(\FF(\phi))\neq 0$
for some $x^{\vee}\in U^{\vee}\cap\bW^{\vee}_0$.
\ed

\section{Global methods}\label{globalsec}

In this section $F$ is a number field, $\A$ is the corresponding
ring of adeles. For every finite set of places $S$ we denote by
$\A^{S}$ (resp. $\A_S$)
the restricted product (resp. the usual product)
of $F_v$ over places $v\not\in S$ (resp. $v\in S$) and by $a\mapsto a^S$
(resp. $a\mapsto a_S$) the corresponding projection from $\A$. 
We fix an algebraic closure $\ov{F}$ of $F$ and
denote by $\Ga$ the Galois group of $\ov{F}$ over $F$.
For every place $v$ of $F$ we denote by $F_v$ the corresponding
local field and by $\Ga_v$ the local Galois group at $v$.
For a reductive group $\bH$ over $F$ we denote by
$\ker^1(F,\bH)$ the preimage of the trivial element under the natural
map
$$H^1(F,\bH)\ra \oplus_v H^1(F_v,\bH).$$
The Hasse principle states that for $\bH$ semisimple and simply connected,
$\ker^1(F,\bH)$ is trivial (see e.g., \cite{PR}).
More generally, for arbitrary
connected reductive group $\bH$, Kottwitz constructed
a bijection
\begin{equation}\label{ker1bij}
\ker^1(F,\bH)\simeq\ker^1(F,Z(\hat{\bH}))^D
\end{equation}
where $Z(\hat{\bH})$ is the centre of the Langlands dual group
(see \cite{K1}, (4.2.2)).
We denote by $\SS(\bV(\A))$ 
the space of Schwartz-Bruhat functions on $\bV(\A)$.
Let $\th$ be the distribution on
$\SS(\bV(\A))$ defined by
$$\th(\phi)=\sum_{x\in\bV(F)}\phi(x),$$
$\th^{\vee}$ be the similar distribution on $\SS(\bV^{\vee}(\A))$.
We set
$$\Th(\phi)=\int_{g\in \bG(\A)/\bG(F)}\th^g(\phi)|dg|$$
where $\th\mapsto\th^g$ denotes the action of $g$ on the distribution $\th$
(provided that the integral converges). Note that
if $\bG$ is anisotropic over $F$ then $\bG(\A)/\bG(F)$ is compact,
so $\Th(\phi)$ is always well-defined in this case.
By Poisson summation formula we have $\FF(\th)=\th^{\vee}$, 
where $\FF$ is the Fourier transform. Hence, for a function
$\phi\in\SS(\bV(\A))$ we have
$$\Th(\FF(\phi))=\Th(\phi)$$
provided that the integral defining $\Th(\phi)$ converges.
In this section we will apply this equality to compare
information about orbital integrals of a function and of its Fourier
transform. In this way we will obtain global proofs of Theorems
\ref{antithm} and \ref{secondmainthm}.

\subsection{Global Kottwitz invariant}

Following Kottwitz we are going
to rewrite the distribution $\Th$ evaluated on sufficiently nice
functions in stably invariant terms.
The main ingredient required for this is the
global Kottwitz invariant defined
in \cite{K1}, \cite{K2}.

Let $\bH$ be a connected reductive group $\bH$ over $F$.
For every place $v$ of $F$ there is a map
$$\inv_v:H^1(F_v,\bH)\ra A(\bH/F_v)=\pi_0(Z(\hat{\bH})^{\Ga_v})^D$$
(see \ref{localKotsec}). Now if we set
$A(\bH/F)=\pi_0(Z(\hat{\bH})^{\Ga})^D$, then for every place $v$
we have a natural homomorphism $r_v:A(\bH/F_v)\ra A(\bH/F)$
induced by the embedding
$Z(\hat{\bH})^{\Ga}\ra Z(\hat{\bH})^{\Ga_v}$.
Thus, we can define a canonical map
$$\inv:\oplus_v H^1(F_v,\bH)\ra A(\bH/F)$$
by setting
$\inv((c_v))=\prod_v r_v\inv_v(c_v).$

The main result about the map $\inv$ is the exactness of the
sequence
$$H^1(F,\bH)\ra \oplus_v H^1(F_v,\bH)\stackrel{\inv}{\ra}
A(\bH/F).$$

\begin{lem}\label{surjKot} Assume that either $\bH$ is semisimple
or $\bH_{F_v}$ is anisotropic. Then the homomorphism
$r_v:A(\bH/F_v)\ra A(\bH/F)$ is surjective.
\end{lem}

\Pf . This follows from the fact that in both cases
$Z(\hat{\bH})^{\Ga_v}$ is finite.
\ed

Now let $\bV$ be a nice representation of $\bG$ defined over $F$.
Then for every $x\in\bV_0(F)$ we can identify 
$$\bG(\A)\backslash\OO_x(\A)$$ with the kernel of the map
$$H^1(\Ga,\bH_x(\ov{\A}))\ra H^1(\Ga,\bG(\ov{\A}))$$
where $\ov{\A}=\A\otimes_F \ov{F}$.
Thus, we can restrict the global Kottwitz invariant
$\inv:H^1(\Ga,\bH_x(\ov{\A}))=\oplus_v H^1(F_v,\bH)\ra A(\bH/F)$
to a map
$$\inv(x,\cdot):\bG(\A)\backslash\OO_x(\A)\ra A(\bH_x/F).$$
We claim that this function takes 
value $0\in A(\bH_x/F)$ precisely on the set of
$\bG(\A)$-orbits of $F$-rational points in $\OO_x$.
Indeed, consider the following commutative diagram
\begin{equation}\label{bigdiagram}
\begin{array}{ccccc}
\bG(F)\backslash\OO_x(F)&\lrar{}& H^1(F,\bH_x)&\lrar{}&H^1(F,\bG)\\
\ldar{}&&\ldar{}&&\ldar{i_{\bG}}\\
\bG(\A)\backslash\OO_x(\A)&\lrar{}&H^1(\Ga,\bH_x(\ov{\A}))&\lrar{}&
H^1(\Ga,\bG(\ov{\A}))\\
&&\ldar{}\\
&&A(\bH_x/F)
\end{array}
\end{equation}
with exact central vertical column. The Hasse principle for $\bG$
(which is simply connected) implies
that the map $i_{\bG}$ is injective. Now our claim follows
by an easy diagram chase.

Let $\bH=\bH_x$ for some $x\in\bV_0(F)$.
Recall that we have defined
in section \ref{localKotsec} a homomorphism
$\sign_v:A(\bH/F_v)\ra \{\pm 1\}$ for every place $v$ such that
$\eps_{\bH_{F_v}}=\sign_v\circ \inv_v$.
We claim that there exists a canonical homomorphism
$$\sign_F:A(\bH/F)\ra\{\pm 1\}$$
such that $\sign_v=\sign_F\circ r_v$ for every place $v$.
Indeed, all the local homomorphisms $\sign_v$ factor through
$A(\bH_{ad}/F_v)$, it suffices to define $\sign_F$ in the case when
$\bH$ is semisimple. This is done in absolutely
the same way as in the local case. Namely, we consider the
exact sequence $1\ra\bC\ra\bH_{sc}\ra\bH\ra 1$ defined over $F$,
where $\bH_{sc}$ is simply connected. The unique morphism of this
exact sequence to $1\ra\{\pm 1\}\ra\wt{\bH}\ra\bH\ra 1$ induces
a homomorphism $\bC\ra\{\pm 1\}$ defined over $F$, that can be considered
as a character $A(\bH/F)\ra\{\pm 1\}$. This is our $\sign_F$.
The compatibility with the local construction is obvious.

\subsection{Stabilization}

\begin{thm} Assume that $\bG$ is a simply connected
semisimple group over $F$,
$\rho:\bG\ra\Aut(\bV)$ is a nice representation defined over $F$.
Assume also that either $\bG$ is anisotropic over $F$, or the
generic stabilizer $\bH$ is semisimple.
Let $f\in\SS(\bV(\A))$ be a function such that
$\supp(f)\cap \bG(\A)\bV(F)\subset\bV_0(\A)$. 
Then
\begin{equation}\label{trfor}
\Th(f)=\sum_{x\in\bG(\ov{F})\backslash\bV_0(F)}
\sum_{\kappa\in A(\bH_x/F)^D}
\int_{a\in\OO_x(\A)}\kappa(\inv(x,a))f(a)|\om_x|,
\end{equation}
where $\bG(\ov{F})\backslash\bV_0(F)$ denotes a set
of representatives of stable $\bG$-equivalence classes of points
in $\bV_0(F)$;
for every $x\in\bV_0(F)$ we choose a $\bG$-invariant
top-degree form $\om_x$ on $\OO_x$ defined over $F$.
\end{thm}

\Pf . For brevity we will omit $|\om_x|$ in the integrals below.
First, we claim that the infinite sum in the right-hand side of
(\ref{trfor}) is absolutely convergent. Indeed, for every point
$x\in\bV_0(F)$ there exists a Zariski neighborhood
$\bU_x\subset \bL_x$ in a linear subspace $\bL_x\subset\bV^{\bH_x}$
such that the map $\bG/\bH_x\times \bU_x\ra\bV_0$ is \'etale.
In particular, the subsets $\bG\bU_x$ form a Zariski open covering of
$\bV_0$, so a finite number of them cover $\bV_0$. Therefore, it
suffices to prove the absolute convergence of
$$\sum_{u\in\bU_x(F)}\int_{a\in\OO_u(\A)}f(a).$$
Now from the fact that the map $\bG/\bH_x\times\bL_x\ra \bV$
is linear in the second argument it is easy to deduce that
the function $u\mapsto \int_{a\in \OO_u(\A)}f(a)$ on $\bU_x(\A)$ is rapidly
decreasing at infinity, which implies our claim.

For every $x\in\bV_0(F)$ we have
$$\sum_{\kappa\in A(\bH_x/F)^D}\int_{a\in\OO_x(\A)}\kappa(\inv(x,a))f(a)=
|A(\bH_x/F)|\cdot \int_{a\in\bG(\A)\OO_x(F)}f(a)$$
since $\inv(x,a)=0$ if and only if $a\in\bG(\A)\OO_x(F)$.
Now we have a finite covering
$$\disj_{y\in\bG(F)\backslash\OO_x(F)}\bG(\A)y\ra\bG(\A)\OO_x(F).$$
We claim that the degree of this covering is equal to
$|\ker^1(F,Z(\hat{\bH}))|$. Indeed, this degree is equal to the number of
$y\in\bG(F)\backslash\OO_x(F)$ such that $\bG(\A)x=\bG(\A)y$.
Such $y$'s correspond to the classes in
$\ker(H^1(F,\bH_x)\ra H^1(F,\bG))$ that have
trivial restriction at every place. Since $\bG$ is simply connected,
by the Hasse principle the elements in $H^1(F,\bG)$ with trivial
restrictions at all places are trivial. Therefore, 
our set of cohomology classes coincides with $\ker^1(F,\bH_x)$.
Using (\ref{ker1bij}) we obtain that
$|\ker^1(F,\bH_x)|=|\ker^1(F,Z(\hat{\bH}_x))|$. It follows that
$$\sum_{\kappa\in A(\bH_x/F)^D}\int_{a\in\OO_x(\A)}\kappa(\inv(x,a))f(a)=
\frac{|A(\bH_x/F)|}{|\ker^1(F,Z(\hat{\bH}_x))|}
\sum_{y\in\bG(F)\backslash\OO_x(F)}\int_{a\in\bG(\A)y}f(a).$$
Now for $y\in\bG(F)$ we have
$$\int_{a\in\bG(\A)y}f(a)=\int_{g\in\bG(\A)/\bH_y(\A)}f(gy)=
\tau(\bH_y)^{-1}\int_{g\in\bG(\A)/\bH_y(F)}f(gy)$$
where $\tau(\bH_y)=\vol(\bH_y(\A)/\bH_y(F))$ is the Tamagawa
number of $\bH_y$.
As Kottwitz showed, the Tamagawa number does not
change if we pass to an inner form. More precisely, we have
$$\tau_{\bH_y}=\frac{|A(\bH_y/F)|}{|\ker^1(F,Z(\hat{\bH_y}))|}$$
(see \cite{K1}, (5.1.1) or the introduction to \cite{K-TN}).
Since $\bH_y$ is an inner form of $\bH_x$, we can replace
$\bH_y$ by $\bH_x$ in the RHS. Hence, we obtain
$$\sum_{\kappa\in A(\bH_x/F)^D}\int_{a\in\OO_x(\A)}\kappa(\inv(x,a))f(a)=
\sum_{y\in\bG(F)\backslash\OO_x(F)}\int_{g\in\bG(\A)/\bH_y(F)}f(gy).$$
Thus, the RHS of the formula
(\ref{trfor}) can be rewritten as follows:
\begin{align*}
&\sum_{x\in\bG(\ov{F})\backslash\bV_0(F)}
\sum_{y\in\bG(F)\backslash\OO_x(F)}\int_{g\in\bG(\A)/\bH_y(F)}f(gy)=
\sum_{y\in\bG(F)\backslash\bV_0(F)}\int_{g\in\bG(\A)/\bH_y(F)}f(gy)=\\
&\int_{g\in \bG(\A)/\bG(F)}\sum_{y\in\bG(F)\backslash\bV_0(F)}
\sum_{x\in\bG(F)y}f(gx)=
\int_{g\in \bG(\A)/\bG(F)}\sum_{x\in\bV_0(F)}f(gx)
\end{align*}
which is precisely the LHS of (\ref{trfor}).
\ed

\subsection{Global proofs}\label{globalpfsec}

We are going to combine
the obtained local and global information about nice representations
with the stabilization formula (\ref{trfor}) to prove theorems \ref{secondmainthm}
and \ref{antithm}.

\noindent
{\it Proof of Theorem \ref{secondmainthm}.}
We will only prove that if $\phi\sim\phi'$ then
$\FF(\phi)\sim_{\eps}\kappa(\bG,\bG')\FF(\phi')$.
The proof of the converse statement is absolutely analogous.
For convenience we divide the proof into several steps.
In step 1 we will extend our local data to the global one
in an appropriate way. In the (crucial) step 2
we apply the stabilization formula to deduce the analogue of
our statement for the product of local fields at two places
of our global field. Finally, in step 3 we will
deduce the statement for one local field.
Let us rename the data $(\bG,\bZ,\bV,\a,\bG',\bV')$ into
$(\bG_E,\bZ_E,\bV_E,\a_E,\bG'_E,\bV'_E)$ to reflect the fact that they
are defined over $E$.

\noindent
STEP 1.
We start by choosing a number field $F$, the data $(\bG,\bZ,\bV)$ defined
over $F$, and a place $v_0$ of $F$ such that
$F_{v_0}=E$, $\bG_{v_0}=\bG_E$, $\bZ_{v_0}=\bZ_{E}$, $\bV_{v_0}=\bV_E$.
In addition we can assume that there are two finite places $v_1$, $v_2$ 
(different from $v_0$) and a real place $v_{\infty}$ of $F$
such that $\bG_{v_1}$ and $\bG_{v_2}$ are split and
$\bG(F_{v_{\infty}})$ is compact.
Note that the fact that $\bV_E$ is a nice representation
implies that the representations
$\bV$ and $\bV_v$ for all places $v$ are nice.

Let us choose a point $x_{v_0}\in\bV_0(F_{v_0})$ such that
the cohomology class $\a_E$ comes from a cohomology class
$\b_E\in H^1(\Ga_{v_0},\bH_{x_{v_0}}/\bZ_{v_0})$ (if no
such point exist then the statement of the theorem is empty).
In the case when the generic stabilizer is commutative
we can also choose a point $x_{v_1}\in\bV_0(F_{v_1})$ such that
$\bH_{x_{v_1}}$ is anisotropic (this is possible by
lemma \ref{anisotrlem}).  We can choose  a global point $x\in\bV_0(F)$
that approximates $x_{v_0}$ and $x_{v_1}$ (resp. $x_{v_0}$ if the
generic stabilizer is non-commutative) well enough,
so that for
$\bH=\bH_x$ we have a class $\b_E\in H^1(\Ga_{v_0},\bH_{v_0}/\bZ_{v_0})$
inducing $\a_E$, and in the case $\bH$ is commutative
we also have that $\bH_{v_1}^0$ is anisotropic. Let us set
$\bK=\bH/\bZ$. Then $\bK_{v_1}$ is either semisimple or an anisotropic
torus, hence, by lemma \ref{surjKot} the homomorphism
$A(\bK/F_{v_1})\ra A(\bK/F)$ is surjective.
Recall that we have an exact sequence
$$H^1(F,\bK)\ra\oplus_v H^1(F_v,\bK)\ra A(\bK/F).$$
By surjectivity of the map $H^1(F_{v_1},\bK)\ra A(\bK/F)$
there exists
an element $\b\in H^1(F,\bK)$ such that $\b_{v_0}=\b_E$ and
$\b_v=0$ for $v\neq v_0, v_1$. Let $\a\in H^1(F,\bG/\bZ)$ be
the class induced by $\b$, and let $(\bG',\bV')$ be the twist of $\bG$
by $\a$. Then $\a_{v_0}=\a_E$, while
the restrictions of $\a$ to all places other than $v_0$ and $v_1$
are trivial. Thus we have $(\bG'_{v_0},\bV'_{v_0})=
(\bG'_E,\bV'_E)$ and $(\bG'_v,\bV'_v)=(\bG_v,\bV_v)$ for
$v\neq v_0, v_1$.
In particular, $\bG'_{v_2}=\bG_{v_2}$ is split over $F_{v_2}$
and $\bG'(F_{v_{\infty}})=\bG(F_{v_{\infty}})$ is compact.
The latter condition implies that both groups $\bG$ and $\bG'$ are
anisotropic over $F$, so the distribution $\Th$ is
defined on all functions in $\SS(\bV(\A))$ (resp. $\SS(\bV'(\A))$,
$\SS(\bV^{\vee}(\A))$ and $\SS((\bV')^{\vee}(\A))$).

\noindent
STEP 2. Set $S=\{v_0,v_1\}$. We are going to show
that for every pair of functions $\phi_S\in\SS(\bV(\A_S))$,
$\phi'_S\in\SS(\bV'(\A_S))$ such that $\phi_S\sim \phi'_S$, one has
$\FF(\phi_S)\sim_{\eps}\FF(\phi'_S)$. For this we have to check that 
for every $x^{\vee}_S\in\bV_0^{\vee}(\A_S)$ one has
$$\delta^{\eps}_{\OO_{x^{\vee}_S},\om}(\FF(\phi_S))=
\delta^{\eps}_{\OO'_{x^{\vee}_S},i_*\om}(\FF(\phi'_S)).$$

Since $\FF(f_S)\in\SS(\bV(\A_S))$ the functions $y_S\mapsto
\delta^{\eps}_{\OO_{y_S},\om}(\FF(\phi_S))$ and
$y_S\mapsto\delta^{\eps}_{\OO'_{y_S},i_*\om}(\FF(\phi'_S))$
on $\bV_0^{\vee}(\A_S)$ are
locally constant (see lemma \ref{nonempty}).
Let $U_S$ be an open neighborhood
of $x^{\vee}_S$ in $\bV_0^{\vee}(\A_S)$ on which these two functions are
constant.

Since $\bG_{v_2}$ is split over $F_{v_2}$, we can
apply the construction of section \ref{localconstrsec} and
proposition \ref{mainloc}, to construct
a stable function $\phi_{v_2}\in\SS(\bV_0(F_{v_2}))$
and a non-empty open subset $U_{v_2}\subset\bV^{\vee}_0(F_{v_2})$
such that the restriction of $I(\FF(\phi_{v_2}))$ to
$\cup_{y\in U_{v_2}}\OO_{y}(F_{v_2})$ is $\eps$-stable
and everywhere non-vanishing.
In particular, for $y\in U_{v_2}$ we have
$\delta^{\eps}_{\OO_{y}}(\FF(\phi_{v_2}))\neq 0$.
If $\bH$ is commutative, then by lemma \ref{anisotrlem}
(and the remark after it) we can in addition assume that all points
in $U_{v_2}$ and in the support of $\phi_{v_2}$ have 
anisotropic stabilizer.

We can find a point $x^{\vee}\in \bV_0^{\vee}(F)$
such that $x^{\vee}\in U_{S}$ and $x^{\vee}\in U_{v_2}$.
Then we have
$$\delta^{\eps}_{(\OO_{x^{\vee}})_{\A_S}}(\FF(\phi_S))=
\delta^{\eps}_{\OO_{x^{\vee}_S}}(\FF(\phi_S)),$$
$$\delta^{\eps}_{(\OO'_{x^{\vee}})_{\A_S}}(\FF(\phi'_S))=
\delta^{\eps}_{\OO'_{x^{\vee}_S}}(\FF(\phi'_S)),$$
$$\delta^{\eps}_{(\OO_{x^{\vee}})_{F_{v_2}}}(\FF(\phi_{v_2}))\neq 0.$$

Let us set $S'=S\cup\{v_2\}$.
We claim that for every 
function $\phi^{S'}\in\SS(\bV(\A^{S'}))$ one has
$$\Theta(\phi_S\otimes\phi_{v_2}\otimes \phi^{S'})=
\Theta(\phi'_S\otimes\phi_{v_2}\otimes \phi^{S'}).$$
Indeed, since $\phi_{v_2}$ is supported on points
in $\bV_0(F_{v_2})$ we can apply
formula (\ref{trfor}) to compute
$\Theta(\phi_S\otimes\phi_{v_2}\otimes \phi^{S'})$.
The RHS of this formula is the sum over $x\in\bV_0(F)$ and over
$\kappa\in A(\bH_x/F)^D$ of terms 
$$\int_{a\in\OO_x(\A)}\kappa(\inv(x,a))
\phi_{S}(a_{S})\phi_{v_2}(a_{v_2})\phi^{S'}(a^{S'})|\om_x|=0.$$
Since $\inv(x,a)$ is the product of local terms, 
this integral is equal to the product of the corresponding
local integrals. We claim that if $\kappa\neq 0$
the local integral at $v_2$ is zero. Indeed,
if $\OO_x(F_{v_2})$ does not intersect the support of $\phi_{v_2}$ this is
clear. Otherwise, by lemma \ref{surjKot} the character
$\kappa\circ r_{v_2}$ of $A(\bH_x/F_{v_2})$ is non-trivial
(here we use the fact that the stabilizers of points in
the support of $\phi_{v_2}$ are either semisimple or anisotropic).
Hence,
$$\int_{a_{v_2}\in\OO_x(F_{v_2})}(\kappa\circ r_{v_2})(a_{v_2})
\phi_{v_2}(a_{v_2})|\om_x|=0$$
since $\ov{\phi}_{v_2}$ is stable.
Thus, the formula (\ref{trfor}) in this case takes form
$$\Theta(\phi_S\otimes\phi_{v_2}\otimes \phi^{S'})=
\sum_{x\in\bG(\ov{F})\backslash\bV_0(F)}
\int_{a\in\OO_x(\A)}\phi_S(a_S)\phi_{v_2}(a_{v_2})
\phi^{S'}(a^{S'})|\om_x|.$$
Similar formula holds for
$\Theta(\phi'_S\otimes\phi_{v_2}\otimes \phi^{S'})$.
Now our claim follows immediately from the condition $\phi_S\sim \phi'_S$.

Applying the Fourier transform, we obtain
\begin{equation}\label{inneraux}
\Theta(\FF(\phi_S)\otimes\FF(\phi_{v_2})\otimes \FF(\phi^{S'}))=
\Theta(\FF(\phi'_S)\otimes\FF(\phi_{v_2})\otimes \FF(\phi^{S'})).
\end{equation}

Let us choose one more finite place
$v_3\not\in S'$ and a function $\phi_{v_3}$
in $\SS(\bV(F_{v_3}))$ 
such that $\FF(\phi_{v_3})$ is supported on
$\bV^{\vee}_0(F_{v_3})$ and
$\delta^{\eps}_{(\OO_{x^{\vee}})_{F_{v_3}}}(\FF(\phi_{v_3}))\neq 0$.
One can define such $\phi_{v_3}$ by setting
$$\FF(\phi_{v_3})=\eps\cdot\delta_{KU_0}$$
where $U_0$ is a small compact neighborhood of $x^{\vee}$ in
$F_{v_3}$-points of a linear slice for $\bG$-action,
$K$ is a non-empty compact
in $\OO_{x^{\vee}}(F_{v_3})$ intersecting $\bG(F_{v_3})$-orbits by the
sets of equal volumes (see section \ref{localconstrsec}
for similar constructions).

Let $S''=S'\cup\{v_3\}$.
Since the function $\FF(\phi_{v_3})$ is supported on
$\bV^{\vee}_0(F_{v_3})$
we can apply formula (\ref{trfor}) to calculate
$\Theta(\FF(\phi))$ and $\Theta(\FF(\phi'))$ where
$$\phi=\phi_S\otimes\phi_{v_2}\otimes\phi_{v_3}\otimes\phi^{S''},$$
$$\phi'=\phi'_S\otimes\phi_{v_2}\otimes\phi_{v_3}\otimes\phi^{S''},$$
for some $\phi^{S''}\in\SS(\bV(\A^{S'}))$.
Now the idea is to choose $\phi^{S''}$ in such a way that 
all the terms in the RHS of (\ref{trfor}) (applied to $\FF(\phi)$ and
$\FF(\phi')$) corresponding to points of
$\bV^{\vee}_0(F)$ which are not stably equivalent to $x^{\vee}$ vanish.
Indeed, let
$C\subset\bV^{\vee}(\A_{S''})$ (resp. $C'\subset\bV^{\vee}(\A_{S''})$)
be the support of $\FF(\phi_S\otimes\phi_{v_2}\otimes \phi_{v_3})$
(resp. of $\FF(\phi'_S\otimes\phi_{v_2}\otimes \phi_{v_3})$).
Set $D=\bV^{\vee}/\bG(F)\cap p_{\A_{S''}}(C\cup C')$,
where $p:\bV^{\vee}\ra\bV^{\vee}/\bG$ is the natural projection
(intersection is taken in $\bV^{\vee}/\bG(\A_{S''})$).
Since $\bV^{\vee}/\bG(F)$ is discrete
in $\bV^{\vee}/\bG(\A)$ and $C\cap C'$ is compact, it follows that $D$ is discrete in
$\bV^{\vee}/\bG(\A^{S''})$.
Note that $p(x^{\vee})$ belongs to $D$. Therefore, we can choose
$\phi^{S''}$ in such a way that its support is disjoint from 
$p_{\A^{S''}}^{-1}(D\setminus p(x^{\vee}))$, while
$\delta^{\epsilon}_{(\OO_{x^{\vee}})_{\A^{S''}}}(\FF(\phi^{S''}))\neq 0$.
Now applying (\ref{trfor}) we get
\begin{align*}
&\Theta(\FF(\phi))=\\
&\sum_{\kappa\in A(\bH_{x^{\vee}}/F)^D}\int_{a\in\OO_{x^{\vee}}(\A)}
\kappa(\inv(x^{\vee},a))\FF(\phi_S)(a_S)\FF(\phi_{v_2})(a_{v_2})
\FF(\phi_{v_3})(a_{v_3})\FF(\phi^{S''})(a^{S''})|\om_x|.
\end{align*}
Since $\bH_{x^{\vee}}$ is either semisimple or anisotropic over $F_{v_2}$,
while the function $I(\FF(\phi_{v_2}))$ is $\eps$-stable
on $\OO_{x^{\vee}}(F_v)$, it follows that
the local integral at $v_2$ vanishes unless 
$\kappa=\eps_F$. The same computation works for $\Theta(\FF(\phi'))$.
Hence, the equality (\ref{inneraux}) for
$\phi^{S'}=\phi_{v_2}\otimes \phi^{S''}$ reduces to
\begin{align*}
&\delta^{\eps}_{(\OO_{x^{\vee}})_{\A_S}}(\FF(\phi_S))\cdot
\delta^{\eps}_{(\OO_{x^{\vee}})_{F_{v_2}}}(\FF(\phi_{v_2}))\cdot
\delta^{\eps}_{(\OO_{x^{\vee}})_{F_{v_3}}}(\FF(\phi_{v_3}))\cdot
\delta^{\eps}_{(\OO_{x^{\vee}})_{\A^{S''}}}(\FF(\phi^{S''}))=\\
&\delta^{\eps}_{(\OO'_{x^{\vee}})_{\A_S}}(\FF(\phi'_S))\cdot
\delta^{\eps}_{(\OO_{x^{\vee}})_{F_{v_2}}}(\FF(\phi_{v_2}))\cdot
\delta^{\eps}_{(\OO_{x^{\vee}})_{F_{v_3}}}(\FF(\phi_{v_3}))\cdot
\delta^{\eps}_{(\OO_{x^{\vee}})_{\A^{S''}}}(\FF(\phi^{S''})).
\end{align*}
Therefore, we get 
$$\delta^{\eps}_{(\OO_{x^{\vee}})_{\A_S}}(\FF(\phi_S))
=\delta^{\eps}_{(\OO'_{x^{\vee}})_{\A_S}}(\FF(\phi'_S))$$
which implies our statement.

\noindent
STEP 3. Applying proposition \ref{innerprop} for the place
$v_1$ we construct functions
$\phi_{v_1}\in\SS(\bV_0(F_{v_1}))$ 
and $\phi'_{v_1}\in\SS(\bV'_0(F_{v_1}))$
such that $\phi_{v_1}\sim \phi'_{v_1}$, while
$$\delta^{\eps}_{\OO,\om}(\FF(\phi_{v_1}))
=\kappa(\bG,\bG')\delta^{\eps}_{\OO',i_*\om}(\FF(\phi'_{v_1}))\neq 0$$
for some orbit $\OO\subset (\bV_0)_{v_1}$. Now let
$\phi_0\in\SS(\bV(F_{v_0}))$, $\phi'_0\in\SS(\bV'(F_{v_0}))$ be
a pair of functions such that $\phi_0\sim \phi'_0$. Then Step 2 applied
to $\phi_0\otimes \phi_{v_1}$ and 
$\phi'_0\otimes \phi'_{v_1}$
implies that $\FF(\phi_0)\sim_{\eps}\kappa(\bG,\bG')\FF(\phi'_0)$.
\ed

\begin{rem} To generalize the above proof to the case when
$\bG$ is not necessarily simple, it would be enough to prove the
following version of lemma \ref{anisotrlem} for $\bG$: if $\bG$ is split
and $\rho$ is nice, then there exists a point $x\in\bV_0(E)$ such that
the connected component of $Z(\bH_x)$ is anisotropic. 
Indeed, the lemmas \ref{maininner} and \ref{surjKot} used in the proof
can be easily generalized to the case when $Z(\bH_x)^0$ is anisotropic.
The rest of the proof does not use the assumption that $\bG$ is simple.
\end{rem}

\noindent
{\it Proof of theorem \ref{antithm}.}
We will only prove the inclusion $\FF(\SS^{st})\subset\SS^{st}_{\eps}$.
The proof of the inverse inclusion is absolutely analogous.
Thus, we have to prove that if $\phi_0$ is a stable function on $\bV(E)$
then for any $x^{\vee}_E\in\bV^{\vee}_0(E)$ the function
$I(\eps\cdot\FF(\phi_0))$ on $\OO_{x^{\vee}_E}(E)$ is constant.
We will split the proof into two steps which are similar to the
first two steps of the previous proof: step 1 consists of
constructing an appropriate global setup, while step 2 is
an application of the stabilization formula and of
theorem \ref{mainthm}.
As before we rename our data $(\bG,\bV)$ into
$(\bG_E,\bV_E)$. 
We denote by $\bH_E\subset\bG_E$ the
stabilizer of $x^{\vee}_E$.

\noindent STEP 1. Let $F$ be a global field, $v_0$ be a place of $F$
such that $F_{v_0}=E$. We want to construct the data $(\bG,\bV,\bH)$
over $F$ such that $\bG_{v_0}=\bG_E$, $\bV_{v_0}=\bV_E$, $\bH_{v_0}$
is $\bG(E)$-conjugate to $\bH_E$ 
and such that in addition the natural homomorphism 
$$A(\bH/E)\ra A(\bH/F)$$
is an isomorphism. 
First, we can find the data $(\bG,\bV)$ over $F$ such that
$\bG_{v_0}=\bG_E$, $\bV_{v_0}=\bV_E$. Now let $x^{\vee}\in\bV^{\vee}(F)$
be a global point sufficiently close to $x^{\vee}_E$, and let $\bH$ be
the stabilizer of $x^{\vee}$. Then
$\bH_{v_0}$ is $\bG(E)$-conjugate to $\bH_E$. Let
$\pi\subset\Aut(Z(\hat{\bH}))$ be the quotient
of the local Galois group $\Ga_{v_0}$, through which it acts on
$Z(\hat{\bH})$. Let us denote by $\Ga'\subset\Ga$ the preimage
of $\pi$ under the natural homomorphism $\Ga\ra\Aut(Z(\hat{\bH}))$,
and let $F'\supset F$ be the finite extension corresponding to
$\Ga'$. Then $\Ga'$ contains $\Ga_v$ (which is
considered as a subgroup in $\Ga$ via some fixed extension of the valuation
$v_0$ to $\ov{F}$). Hence, there is an extension of $v_0$ to a place
$v'_0$ of $F'$ such that $F'_{v'_0}=E$. Furthermore, by construction
we have $A(\bH/E)=A(\bH/F')$. It remains to replace $(F,v_0)$ by
$(F',v'_0)$.

\noindent STEP 2. This step is very similar to the step 2 in the
previous proof.
Let $\kappa_E$ be a non-trivial character of $A(\bH/E)$. We have
to prove that 
$\delta^{\kappa_E\sign}_{\OO_{x^{\vee}_E}}(\FF(\phi_0))=0$.
Let us denote $\bW^{\vee}=(\bV^{\vee})^{\bH}$. Since the function
$y\mapsto \delta^{\kappa_E\sign}_{\OO_y}(\FF(\phi_0))$
of $y\in\bW^{\vee}(E)$ is locally constant, we can choose
a neighborhood $U_{v_0}$ of $x^{\vee}_E$ on which this function is
constant. By our assumption $\kappa_E$ is induced
by some character $\kappa$ of $A(\bH/F)$. For every place $v$ we
denote by $\kappa_v$ the induced character of
$A(\bH/F_{v_1})$.

Let us choose a finite place
$v_1$ of $F$ (different from $v_0$). By proposition
\ref{mainloc} (ii) for a function of the form
$\phi_{v_1}=\delta^{\kappa_{v_1}}_{U(t)}\in\SS(\bV_0(F_{v_1}))$
where $t\in (E^*)^2$ is large enough, one has
$$\delta^{\kappa_{v_1}\sign}_{\OO_{y^{\vee}}}(\FF(\phi_{v_1}))
\neq 0$$ for $y^{\vee}\in U_{v_1}$, where $U_{v_1}\subset
\bW^{\vee}_0(F_{v_1})$ is a non-empty open subset.
Also since the character $\kappa_{v_1}$ is non-trivial by lemma
\ref{surjKot}, the function $\phi_{v_1}$ is antistable.

Let us choose a global point $x^{\vee}\in\bW^{\vee}_0(F)$ such that
$x^{\vee}\in U_{v_0}$ and $x^{\vee}\in U_{v_1}$. Then we have
$$\delta^{\kappa_E\sign}_{\OO_{x^{\vee}_E}}(\FF(\phi_0))=
\delta^{\kappa_E\sign}_{(\OO_{x^{\vee}})_{F_{v_0}}}(\FF(\phi_0)),$$
$$\delta^{\kappa_{v_1}\sign}_{(\OO_{x^{\vee}})_{F_{v_1}}}(\FF(\phi_{v_1}))
\neq 0.$$

Let $v_2$ be one more finite place of $F$.
We can construct a function $\phi_{v_2}\in\SS(\bV(F_{v_2}))$
such that $\FF(\phi_{v_2})$ has support in $\bV_0^{\vee}(F_{v_2})$
and such that
$$\delta^{\kappa'}_{(\OO_{x^{\vee}})_{F_{v_2}}}(\FF(\phi_{v_2}))=0$$
for $\kappa'\neq\kappa_{v_2}\sign$, while 
$$\delta^{\kappa_{v_2}\sign}_{(\OO_{x^{\vee}})_{F_{v_2}}}(\FF(\phi_{v_2}))
\neq 0.$$
Indeed, it suffices to take $\FF(\phi_{v_2})$ to be the function
of the form $\delta^{\kappa_{v_2}\sign}_{U(t)}$ as in section
\ref{localconstrsec}.

Let us denote $S=\{v_0,v_1,v_2\}$. 
We claim that for every 
function $\phi^{S}\in\SS(\bV(\A^{S}))$ one has
$$\Theta(\phi_0\otimes\phi_{v_1}\otimes\phi_{v_2}\otimes \phi^{S})=0.$$
Indeed, since $\phi_{v_1}$ has support in
$\bV_0(F_{v_1})$, we can apply
formula (\ref{trfor}).
Now let $x\in\bV_0(F)$ and let $\kappa'$ be a character of $A(\bH_x/F)$.
If $\kappa'$ is non-trivial then by lemma \ref{surjKot} the induced
character $\kappa'_{v_0}$ of $A(\bH_x/E)$ is non-trivial. Since $\phi_0$
is stable, we get
$$\int_{a\in\OO_x(E)}\kappa'(\inv(x,a))\phi_0(a)|\om_x|=0,$$
so the corresponding term in the RHS of (\ref{trfor}) vanishes.
On the other hand, if $\kappa'=1$ then $\kappa'_{v_1}=1$,
so the corresponding term vanishes by antistability of $\phi_{v_1}$.

Applying the Fourier transform, we obtain
$$\Theta(\FF(\phi_0)\otimes\FF(\phi_{v_1})\otimes \FF(\phi_{v_2})\otimes
\FF(\phi^S))=0.$$
Furthermore, since $\FF(\phi_{v_2})$ has support in $\bV_0^{\vee}(F_{v_2})$
we can apply formula (\ref{trfor}) again.
As in the proof of theorem \ref{mainthm} we can
choose $\phi^{S}$ in such a way that
all the terms in the RHS of (\ref{trfor}) corresponding to points of
$\bV^{\vee}_0(F)$ which are not stably conjugate to $x^{\vee}$ vanish
while
$$\delta_{(\OO_{x^{\vee}})_{\A^S}}(\FF(\phi^S))\neq 0.$$
Now applying (\ref{trfor}) we get
\begin{align*}
&0=\Theta(\FF(\phi_0)\otimes\FF(\phi_{v_1})\otimes\FF(\phi_{v_2})\otimes
\FF(\phi^S))=\\
&\sum_{\kappa'\in A(\bH_{x^{\vee}}/F)^D}\int_{a\in\OO_{x^{\vee}}(\A)}
\kappa'(\inv(x^{\vee},a))\FF(\phi_0)(a_{v_0})\FF(\phi_{v_1})(a_{v_1})
\FF(\phi_{v_2})(a_{v_2})\FF(\phi^{S})(a^{S})|\om_x|.
\end{align*}
By our choice of $\phi_{v_2}$
the local integral at $v_2$ vanishes unless
$\kappa'_{v_2}=\kappa_{v_2}\sign$. By lemma \ref{surjKot}
this condition is equivalent to $\kappa'=\kappa\sign_F$.
Hence, we obtain
that $\delta^{\kappa_E\sign}_{(\OO_{x^{\vee}})_{F_{v_0}}}(\FF(\phi_0))=0$
as required.
\ed

\section{Sign function for the space of symmetric matrices}
\label{prehom}

In this section we will compute the sign function $\eps(\cdot,\cdot)$
for the space of symmetric $n\times n$ matrices $\Sym_n$
over a local field $E$ considered as a representation of $\SL_n$,
where $g\in\SL_n$ acts on $\Sym_n$ by $X\mapsto gXg^t$.
As an application we will derive the formula (\ref{discreq}) in the
case of odd $n$.

Let us denote by $\Sym'_n\subset\Sym_n$ the complement to the
hypersurface $(\det=0)$.
By definition the sign $\eps(A,A')$, where $A,A'\in\Sym'_n(E)$,
is defined when $A$ and $A'$ belong to one $\SL_n(\ov{E})$-orbit,
i.e., when $\det(A)=\det(A')$. 

\begin{prop}\label{signprop}
For a pair of symmetric $n\times n$ matrices $A,A'$
with $\det(A)=\det(A')\neq 0$ one has
$$\eps(A,A')=\eps(q_A,q_{A'})^n$$
where $q_A$ denotes the quadratic form with matrix $A$.
\end{prop}

\Pf . By definition $\eps(A,A')$ is the product of Hasse-Witt invariants
of quadratic forms obtained by restricting the Killing form $Q$ on
$\ssl_n$ to stabilizer subalgebras $\hh_A$ and $\hh_{A'}$ of $A$ and $A'$.
Since $A$ and $A'$ belong to one $\SL_n(\ov{E})$-orbit, the forms
$Q|_{\hh_A}$ and $Q|_{\hh_{A'}}$ have the same determinant modulo squares.
Therefore, their relative Hasse-Witt invariant will not change if
we replace $Q$ by its scalar multiple. Thus, we can do calculation
with $Q(X)=\frac{1}{2}Tr(X^2)$. We claim that
$$\eps(Q|_{\hh_A})=c(\det A)\cdot\eps(q_A)^n$$
where $c$ is a sign depending only on $\det A$ modulo squares.
To prove this formula we notice that both sides do not change
if we replace $A$ by $gAg^t$ where $g\in\GL_n(E)$. Thus, we
can assume that $A$ is diagonal. Let $(a_1,\ldots,a_n)$ be
diagonal entries of $A$. The subalgebra $\hh_A\subset\ssl_n$
consists of matrices $X$ such that $XA+AX^t=0$. Thus, if $X=(x_{ij})$
then we should have $x_{ii}=0$ while $x_{ji}=-\frac{a_j}{a_i}x_{ij}$.
Thus, the quadratic form $Q|_{\hh_A}$ has diagonal matrix in
the natural basis on $\hh_A$ and we have
$$\eps(Q|_{\hh_A})=
\prod_{i<j,k<l:(i,j)<(k,l)}(-\frac{a_j}{a_i},-\frac{a_l}{a_k})$$
where $(i,j)<(k,l)$ denotes the lexicographical order.
A straightforward calculation shows that the RHS is equal to
$$\prod_{i<j}(a_i,a_j)^n\cdot\prod_i (a_i,-1)^{p(n)}\cdot
(-1,-1)^{q(n)}$$
where $p(n),q(n)$ are some polynomials in $n$.
Thus, we get
$$\eps(Q|_{\hh_A})=\eps(q_A)^n\cdot (\det A,-1)^{p(n)}\cdot (-1,-1)^{q(n)}$$
as required.
\ed

Now we can derive the formula (\ref{discreq}) in the $p$-adic case.
It is well-known (see \cite{S},\cite{SS},\cite{Ig}) 
that one has an equation of the form
$$
\FF(\chi(\det))=\tau_{\chi}\cdot
(|\cdot|^{-\frac{n+1}{2}}\chi^{-1})(\det)
$$
where $\tau_{\chi}$ is some $(\GL_n/\{\pm 1\})(E)$-invariant
function on $\Sym'_n(E)$.
The stabilizer subgroup $\bH\subset (\GL_n/\{\pm 1\})$
of a point in $\Sym'_n(E)$ is the group $\O_n/\{\pm 1\}\simeq\SO_n$
(here we use the fact that $n$ is odd). Therefore,
the set of $(\GL_n/\{\pm 1\})(E)$-orbits on $\Sym'_n(E)$ can be
identified with $\ker(H^1(E,\SO_n)\ra H^1(E,\GL_n/\{\pm 1\}))$.
Note that the homomorphism $\SO_n\ra\GL_n/\{\pm 1\}$ factors through
$\SL_n$ (since $n$ is odd). Therefore, the map
$H^1(E,\SO_n)\ra H^1(E,\GL_n/\{\pm 1\})$ is trivial, so there are
two $(\GL_n/\{\pm 1\})(E)$ orbits on $\Sym'_n(E)$ corresponding
to two distinct elements of $H^1(E,\SO_n)$.
These two orbits intersect the subset $\det A=1$ by
two distinct $\SL_n(E)$-orbits. Therefore, the function $\tau_{\chi}$
is determined by its restriction to the subset $\det A=1$.
It remains to notice that the distribution $\chi(\det)$ is stable
for the action of $\SL_n$ on $\Sym_n$. Therefore, by theorem
\ref{mainthm} its Fourier transform is $\eps$-stable. Applying
proposition \ref{signprop} we conclude that
$$\tau_{\chi}(A)=c(\chi)\cdot \eps(q_A)\cdot f(\det A)$$
where $f$ is a function 
of $\det A\mod (E^*)^2$ (unique up to a constant)
such that $A\mapsto\eps(q_A)\cdot f(\det A)$ is
$(\GL_n/\{\pm 1\})(E)$-invariant. It is easy to see that
$(\GL_n/\{\pm 1\})(E)$-orbits on $\Sym'_n(E)$ coincide with
$\GL_n(E)\times E^*$-orbits on it, where $E^*$ acts by rescaling.
Now for $t\in E^*$ we have
$$\eps(t q_A)=(t,t)^{\frac{n(n-1)}{2}}(t,\det A)^{\frac{n-1}{2}}\eps(q_A)=
(t,-1)^{\frac{n(n-1)}{2}}\eps(q_A).$$
It follows that the function
$A\mapsto \eps(q_A)(\det A,-1)^{\frac{n-1}{2}}$ is
$(\GL_n/\{\pm 1\})(E)$-invariant, so we should have
$$\tau_{\chi}(A)=c(\chi)\cdot \eps(q_A)\cdot (\det A,-1)^{\frac{n-1}{2}}$$
which is equivalent to the equation (\ref{discreq}).
The explicit value of the constant $c(\chi)$ can be found in \cite{Sw}.
Note that it can also be determined using the stationary phase
approximation as in section \ref{statsec}.

In conclusion let us show that for $E=\R$ the formula (\ref{discreq}) is still true.
Indeed, we have either $\chi(x)=|x|^s$ or
$\chi(x)=\sgn(x)|x|^s$, so our formula is equivalent to the set of two equalities 
\begin{equation}\label{real1}
\int_{A\in\Sym'_n(\R)}|\det(A)|^s\hat{f}(A)dA=
c_1(s)\int_{B\in\Sym'_n(\R)}(-1)^{\frac{i_B(n-i_B)}{2}}\cdot
|\det(B)|^{-s-\frac{n+1}{2}}f(B)dB,
\end{equation}
\begin{equation}\label{real2}
\begin{array}{l}
\int_{A\in\Sym'_n(\R)}\sgn(\det(A))|\det(A)|^s\hat{f}(A)dA=\\
c_2(s)\int_{B\in\Sym'_n(\R)}(-1)^{\frac{i_B(n-i_B)}{2}}\cdot\sgn(\det(B))
|\det(B)|^{-s-\frac{n+1}{2}}f(B)dB,
\end{array}
\end{equation}
where $i_B$ is the number of negative eigenvalues of $B$,
$f$ is a function from the Schwartz space of $\Sym_n(\R)$,
$\hat{f}$ is its Fourier transform, $c_1$ and $c_2$ are some meromorphic
functions of $s$. These formulas can be deduced from the explicit form
of $\Ga$-matrix for $\Sym_n(\R)$ computed by T.~Shintani in Lemma 15 of
\cite{Sh}. Indeed, let us recall the result of Shintani's computation
in the form convenient for us. Let us denote by $V_i$ the connected component
of $\Sym'_n(\R)$ consisting of matrices with exactly $i$ positive
eigenvalues. Let us denote $\Phi_i(f,s)=\int_{A\in V_i}|\det(A)|^s f(A)dA$.
Then one has the following equation:
$$\Phi_i(\hat{f},s)=c(s)\cdot\sum_{j=0}^n
v_{ij}(s)\Phi_j(f,-s-\frac{n+1}{2})$$
where
$$v_{ij}(s)=\sum_{(\eps_1,\ldots,\eps_n)}
\exp(\frac{\pi \sqrt{-1}}{2}[\sum_{k=1}^j(k+s)\eps_k-
\sum_{k=j+1}^n(k-j+s)\eps_k])$$
where the summation is taken over all $n$-tuples
$(\eps_1,\ldots,\eps_n)=(\pm 1,\ldots,\pm 1)$ such that exactly
$i$ of the $\eps$'s are $+1$. Now we have
$$c_j:=
\sum_{i=0}^n v_{ij}=2^n\cdot\prod_{k=1}^j\cos(\frac{\pi}{2}(k+s))\cdot
\prod_{k=1}^{n-j}\cos(\frac{\pi}{2}(k+s)).$$
Hence, for odd $n$ we have
$$\frac{c_j}{c_{j-1}}=\frac{\cos(\frac{\pi}{2}(j+s))}
{\cos(\frac{\pi}{2}(n+1-j+s))}=(-1)^{\frac{n+1}{2}+j}.$$
Therefore, the vector $(c_j)_{j=0,\ldots,n}$ is proportional
to $((-1)^{j(n-j)})$ which is equivalent to the equation (\ref{real1}).
Similarly, we have
\begin{align*}
&c'_j:=\sum_{i=0}^n (-1)^{n-i} v_{ij}=
\sum_{(\eps_1,\ldots,\eps_n)\in\{\pm 1\}^n}
\prod_{k=1}^n\eps_k\cdot
\exp(\frac{\pi\sqrt{-1}}{2}[\sum_{k=1}^j(k+s)\eps_k-
\sum_{k=j+1}^n(k-j+s)\eps_k])=\\
&(2\sqrt{-1})^n\cdot(-1)^{n-j}\cdot
\prod_{k=1}^j\sin(\frac{\pi}{2}(k+s))\cdot
\prod_{k=1}^{n-j}\sin(\frac{\pi}{2}(k+s)).
\end{align*}
Hence, for odd $n$ the vector $(c'_j)_{j=0,\ldots,n}$ is proportional
to $((-1)^{(j+1)(n-j)})$ which is equivalent to (\ref{real2}).

\vspace{3mm}

{\sc Department of Mathematics, Harvard University, Cambridge MA 02138

Department of Mathematics and Statistics, Boston University, Boston MA 02215}

{\it E-mail addresses:} kazhdan@@math.harvard.edu, apolish@@math.bu.edu

\end{document}